\documentclass[12pt]{article}

\usepackage{graphicx}

\usepackage{amsmath, amsfonts, amssymb}
\usepackage{epsf}
\newtheorem{theorem}{Theorem}[section]

\newtheorem{proposition}[theorem]{Proposition}
\newtheorem{corollary}[theorem]{Corollary}

\newtheorem{definition}[theorem]{Definition}
\newtheorem{conjecture}[theorem]{Conjecture}
\newtheorem{Problem}[theorem]{Problem}
\newtheorem{question}[theorem]{Question}
\newcommand{\BB}{{\mathbb B}}
\newcommand{\RR}{{\mathbb R}}
\newcommand{\CC}{{\mathbb C}}

\newcommand{\aut}{\hbox{\rm Aut }}

\newcommand{\calC}{{\mathcal C}}

\newfam\msbfam
\font\tenmsb=msbm10 scaled \magstep1 \textfont\msbfam=\tenmsb
\font\sevenmsb=msbm7 scaled \magstep1 \scriptfont\msbfam=\sevenmsb
\font\fivemsb=msbm5 scaled \magstep1
\scriptscriptfont\msbfam=\fivemsb
\def\Bbb{\fam\msbfam \tenmsb}

\def\dbar{\overline{\partial}}

\def\ra{\rightarrow}

\def\ss{\subseteq}

\def\a{\alpha}
\def\Re{\hbox{\rm Re}\,}
\def\Im{\hbox{\rm Im}\,}
\def\Aut{\hbox{Aut}\,}

\def\O{\Omega}
\def\bO{\partial \Omega}
 \def\HollowBox #1#2{{\dimen0=#1 \advance\dimen0 by -#2
   \dimen1=#1 \advance\dimen1 by #2
    \vrule height #1 depth #2 width #2
    \vrule height 0pt depth #2 width #1
        \llap{\vrule height #1 depth -\dimen0 width \dimen1}%
      \hskip -#2
      \vrule height #1 depth #2 width #2}}
 \def\BoxOpTwo{\mathord{\HollowBox{6pt}{.4pt}}\;}
\def\endpf{\hfill $\BoxOpTwo$}

\newcommand{\set}[1]{\{#1\}}
\newcommand{\im}{{\rm Im}\,}
\newcommand{\re}{{\rm Re}\,}

\input epsf

\begin{document}

\title{Complex Scaling and Geometric Analysis \\
of Several Variables}

\author{Kang-Tae Kim \and Steven G. Krantz}

\maketitle
\setcounter{section}{-1}

\section{Preliminary Remarks} %
It is a classical fact that there is no Riemann mapping theorem in
the function theory of several complex variables. Indeed, H.
Poincar\'{e} proved in 1906 that the unit ball $B = \{z = (z_1,
z_2) \in \CC^2: |z| \equiv |z_1|^2 + |z_2|^2 < 1\}$ and the unit
bidisc $D^2 = \{z = (z_1, z_2) \in \CC^2: |z_1| < 1, |z_2| < 1\}$
are biholomorphically inequivalent. More recently,
Burns/Shnider/Wells \cite{BSW} and Greene/Krantz \cite{GRK1} have
shown that two smoothly bounded, strongly pseudoconvex domains
(definitions to be discussed below) are generically
biholomorphically inequivalent. In particular, if one concentrates
attention on smoothly bounded domains that are near the unit ball
in some reasonable metric then, with probability one, two randomly
selected domains will be biholomorphically inequivalent. Thus one
seeks substitutes for the Riemann mapping theorem. In particular,
one seeks to classify domains in terms of geometric invariants.
Work of Fefferman \cite{FEF}, Bell \cite{BEL}, Bell/Ligocka
\cite{BELI}, and others has shown us that a biholomorphic mapping
of a reasonable class of smoothly bounded domains will extend to
smooth diffeomorphism of the closures. [More recently, McNeal
\cite{MCN} has announced that a biholomorphic mapping of {\it any}
smoothly bounded domains will continue to a diffeomorphism of the
closures.]

Thus it is possible, at least in principle, to carry out
Poincar\'{e}'s original program of determining differential
biholomorphic invariants on the boundary. Chern and Moser
\cite{CHM} did the initial work in this direction. (See also
\cite{TAN}.) More recent progress has been made by Webster
\cite{WEB1}, Moser \cite{MOS1} and \cite{MOS2},
Moser and Webster \cite{MOW}, Isaev/Kruzhilin
\cite{ISK}, and Ejov/Isaev \cite{EJI}. Another direction, also
inspired by Poincar\'{e}'s work, is to study the automorphism
group of a domain. This is a natural biholomorphic invariant, and
reflects the Levi and Bergman geometry of the domain in a variety
of subtle and useful ways. The purpose of this paper to is develop
some techniques connected with the study of automorphism groups of
bounded domains in $\CC^{n+1}$. In particular, we wish to focus
attention on a powerful technique that has become central in the
subject. This is the method of {\it scaling}. A special case of a
general technique in differential geometry known as {\it
flattening}, scaling is a method for localizing analysis near a
boundary point. This method has been used with considerable
effectiveness to study not only automorphism groups (\cite{PIN},
\cite{FRA}, \cite{GRK2}) but also canonical kernels (\cite{NRSW},
\cite{KRA2}) and other aspects of classical function theory. It is
a far-reaching methodology that has potential applications in many
parts of mathematics. We exhibit in this article several different
contexts and applications in which the scaling point of view is
useful. For our purposes, the theory of automorphism groups is a
convenient venue in which to showcase the scaling technique. But
it should be of interest to mathematicians with many diverse
interests.

\section{Introduction}

By a {\it domain}, we mean a connected open subset in $\CC^{n+1}$
for some positive integer\footnote{It is notationally convenient
in this study to treat domains in $\CC^{n+1}$ rather than
$\CC^n$.} $n$. Throughout this paper, we shall use $z=(z_0, z_1,
\dots, z_n)$ for the coordinates of a point in $\CC^{n+1}$.

Let $\Omega$ be a domain. If its boundary $\partial\Omega$ is a
regularly imbedded $\mathcal C^k$-hypersurface ($k \ge 2$), then
there exists a $\mathcal C^k$ smooth function $\rho:\CC^{n+1} \to
\RR$ such that
\begin{itemize} \item[{\bf (i)}] $ \Omega = \set{z \in
\CC^{n+1} : \rho(z)<0}$
\item[{\bf (ii)}] $\nabla \rho (p) \not= 0$ whenever $p \in
\partial\Omega$. \end{itemize} In such a case, $\Omega$ is
called a domain with $\mathcal C^k$ smooth boundary. In turn,
$\rho$ is called a $\mathcal C^k$ smooth {\it defining
function} for $\Omega$. \medskip Now fix a domain $\Omega \ss
\CC^{n+1}$ with $\mathcal C^2$ boundary and defining function $\rho$.
Let $p \in \partial
\Omega$. We say that $w \in \CC^{n+1}$ is a {\it complex tangent
vector} to $\partial \Omega$ at $p$ if
$$
\sum_{j=0}^{n} \frac{\partial \rho}{\partial z_j}(p) w_j = 0 \, .
$$
We write $w \in {\mathcal T}_p(\partial \Omega)$. [Observe that
$w$ is a {\it real tangent vector}, which is the classical notion
of tangency from differential geometry, if $\hbox{Re}\,
\sum_{j=0}^{n} \frac{\partial \rho}{\partial z_j}(p) w_j = 0$. We
write in this case $w \in T_p(\partial \Omega)$.] The {\it complex
normal directions} are the directions in $T_p(\partial \Omega)$
which are complementary to ${\mathcal T}_P(\bO)$.

We say that $\partial \Omega$ is (weakly) Levi pseudoconvex at
$p$ if
$$
\sum_{j, k = 0}^{n} \frac{\partial^2 \rho}{\partial z_j \partial
\overline{z}_k} (p) w_j \overline{w}_k \geq 0 \eqno (*)
$$
for every complex tangent vector $w$ at $p$. The point $p$ is {\it
strictly} or {\it strongly} pseudoconvex if the inequality in
$(*)$ is strict whenever $0 \ne w \in {\mathcal T}_p(\partial
\Omega)$. If each point of $\partial \Omega$ is pseudoconvex then
the domain is said to be pseudoconvex; if each point of $\partial
\Omega$ is strongly pseudoconvex then the domain is said to be
strongly pseudoconvex. We note in passing that there is a more
general notion of pseudoconvexity due to Hartogs, and which
utilizes the theory of plurisubharmonic functions. We shall have
no use for that concept here, but see \cite{KRA1}. It is worth
noting (see \cite{KRA1}) that if $p \in \partial \Omega$ is a
point of strong pseudoconvexity then there is a choice of defining
function $\widetilde{\rho}$ so that
$$
\sum_{j, k = 1}^{n+1} \frac{\partial^2 \widetilde{\rho}}{\partial
z_j \partial \overline{z}_k} (p) w_j \overline{w}_k > 0
$$
for {\it every} $w \in \CC^{n+1} \setminus \{0\}$ (not just the
complex tangential $w$). In fact one may go further. There is a
biholomorphic change of coordinates in a neighborhood of $p$ so
that the boundary near $p$ is strongly {\it convex}. This means
that, identifying $z_j = t_{2j-1} + i t_{2j}$, and choosing an
appropriate defining function $\widetilde{\widetilde{\rho}}$, we
have
$$
\sum_{j,k = -1}^{2n} \frac{\partial^2
\widetilde{\widetilde{\rho}}}{\partial t_j \partial t_k} (p) a_j
a_k > 0
$$
for every non-zero real vector $a = (a_1, a_2, \dots, a_{2n})$.
Again see \cite{KRA1} for the details.

\section{The Lore of Automorphism Groups}

Let $\Omega \subseteq \CC^{n+1}$ be a domain. The {\it
automorphism group} of $\Omega$ is the collection of one-to-one,
onto holomorphic mappings $\varphi: \Omega \rightarrow \Omega$. It
is known that the inverse $\varphi^{-1}$ is automatically
holomorphic. See \cite{KRA1}. Such a mapping is also called a {\it
biholomorphic self-map} of $\Omega$. With the binary operation of
composition of mappings, the collection of automorphisms forms a
group. We denote this group by $\hbox{Aut}(\Omega)$. We equip the
automorphism group with the topology of uniform convergence on
compact sets, equivalently the compact-open topology. If we
restrict attention to bounded domains---and in this paper we, for
the most part, do just that---then the group $\Aut(\Omega)$ is in
fact a real Lie group (this follows from work of H. Cartan---see
\cite{KOB1}). It is never a complex Lie group unless it is
discrete. Of special interest are domains having a large or robust
group of automorphisms. It is known (see \cite{BSW}, \cite{GRK1})
that strongly pseudoconvex domains which are rigid---i.e., which
have no automorphisms except the identity---are generic. On the
other hand, every compact Lie group arises as the automorphism
group of some bounded, strongly pseudoconvex domain with real
analytic boundary (see \cite{BDA}, \cite{SAZA}, \cite{GRK4} and
\cite{WIN}). It seems natural, for example, to study a domain with
{\it transitive} automorphism group---this is a domain $\Omega$
with the property that if $P, Q \in \Omega$ are arbitrary then
there is an automorphism $\varphi$ of $\Omega$ such that
$\varphi(P) = Q$. It turns out that the list of such domains is
rather restrictive; our knowledge of such domains is essentially
complete (see \cite{HEL}, \cite{HUC}).\footnote{In fact the only
strongly pseudoconvex domain with transitive automorphism group is
the ball.  This remarkable fact will be discussed below.} Perhaps
geometrically more natural is to consider domains with {\it
noncompact} automorphism group. A very natural and useful
characterization of such domains is contained in the following
classical result of Cartan (for which see \cite{NAR}):

\begin{proposition} \sl
Let $\Omega \ss \CC^{n+1}$ be a bounded domain with noncompact
automorphism group. Then there are a point $p \in \partial
\Omega$, a point $q \in \Omega$, and automorphisms $\varphi_j \in
\Aut(\Omega)$ such that $\varphi_j(q) \ra p$ as $j \ra \infty$.
See Figure 1.
\end{proposition} %

\begin{figure}
\centering
\includegraphics[height=3.25in, width=2.25in]{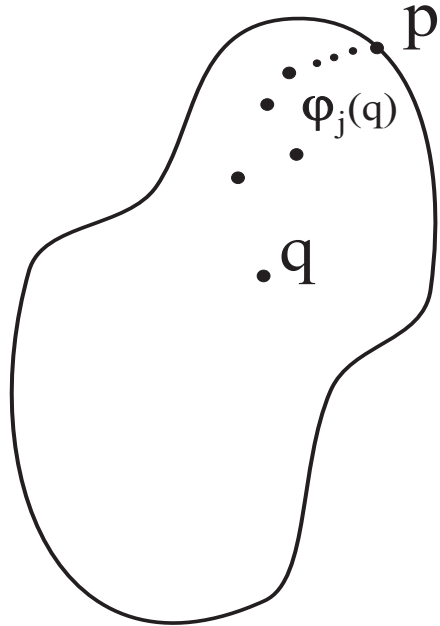}
\caption{Cartan's lemma about noncompact automorphism groups.}
\end{figure}

The point $p$ in the proposition is called a {\it boundary orbit
accumulation point}. It is known, in a variety of concrete senses,
that the Levi geometry of a boundary orbit accumulation point says
a great deal about the domain itself. And conversely. See
\cite{GRK5}, \cite{GRK6}. It is a matter of great interest today
to classify all possible boundary orbit accumulation points. An
important focus of our studies will be (bounded) domains in
$\CC^{n+1}$ with noncompact automorphism group and their boundary
orbit accumulation points. Much is known today about automorphism
groups of domains. In classical studies, mathematicians calculated
the automorphism groups of very particular domains rather
explicitly (see, for instance, the discussion in \cite{KRA1} as
well as \cite{HUA}). Today we have more powerful machinery (the
$\overline{\partial}$-Neumann problem, sheaf theory, Levi
geometry, Lie group methods, techniques of Riemannian geometry,
K\"{a}hler theory) that allow us to make qualitative studies of
broad classes of domains. The papers \cite{GRK3} and \cite{IKR}
provide a broad overview of the types of results that can be
proved with modern techniques. The present paper will introduce
the reader to some of the main themes in the subject.

\section{The Dilatation and Scaling Sequences}

At this stage, we shall only consider the case when the boundary
$\partial\Omega$ is $\mathcal C^k$ smooth, with $k \ge 2$, in an
open neighborhood of $p \in \partial \Omega$. Let $q^\nu =
(q_0^\nu, \ldots, q_n^\nu)$ be a sequence of points in the closure
$\overline{\Omega}$ of the domain $\Omega$ that converges to the
boundary point $p=(p_0, \ldots, p_n)$.

\subsection{Pinchuk's Dilatation Sequence}

We refer to the source \cite{PIN} for basic ideas about Pinchuk
scaling. Applying a holomorphic coordinate change and the implicit
function theorem at $p$, we may assume that $p$ is the origin and
the domain $\Omega$ is represented in an open neighborhood of the
origin by (i.e., has a defining function given by) the defining
inequality
$$
\re z_0 \ge \psi(\im z_0, z_1, \ldots, z_n),
$$
where:
\begin{itemize}
\item[(\romannumeral1)] $\psi \in \mathcal C^k$,
\item[(\romannumeral2)] $\psi(0, \ldots, 0) = 0$, and
\item[(\romannumeral3)] $\nabla \psi|_{(0,\ldots,0)} = (0, \ldots,
0)$.
\end{itemize}

We take $(-1,0, \dots, 0)$ to be the unit outward normal vector at
$p$. Now choose the boundary points $p^\nu = (p_0^\nu, \ldots,
p_n^\nu)$ satisfying
\begin{itemize} \item[(\romannumeral1)] $p_j^\nu = \ q_j^\nu$
for every $j = 1, \ldots, n$, and
\item[(\romannumeral2)] $q_0^\nu - p_0^\nu > 0$
\end{itemize}
for every $\nu = 1,2, \ldots$. Observe that each $p^\nu$ is
uniquely determined by these conditions. See Figure 2.

\begin{figure}
\centering
\includegraphics[height=2.65in, width=2.75in]{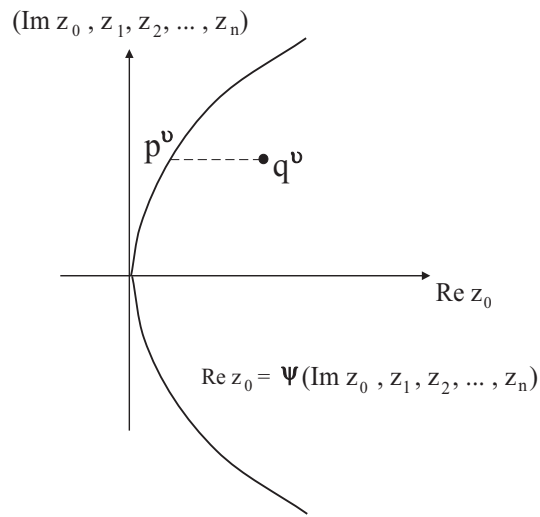}
\caption{A scaling sequence.}
\end{figure}

Then consider the map $A_\nu : \CC^{n+1} \to \CC^{n+1}$ defined by
$\zeta = A_\nu (z)$ in local coordinates with the explicit
expression
\begin{eqnarray*}
\zeta_0  & = & \alpha_0^\nu (z_0 - p_0^\nu) - \sum_{j=1}^{n+1}
\alpha_j^\nu (z_j - p_j^\nu) \\
\zeta_1  & = & z_1 - p_1^\nu \\
 & \vdots &   \\
\zeta_n  & =  & z_n - p_n^\nu \, .
\end{eqnarray*}
See Figure 3.

\begin{figure}
\centering
\includegraphics[height=2.65in, width=2.75in]{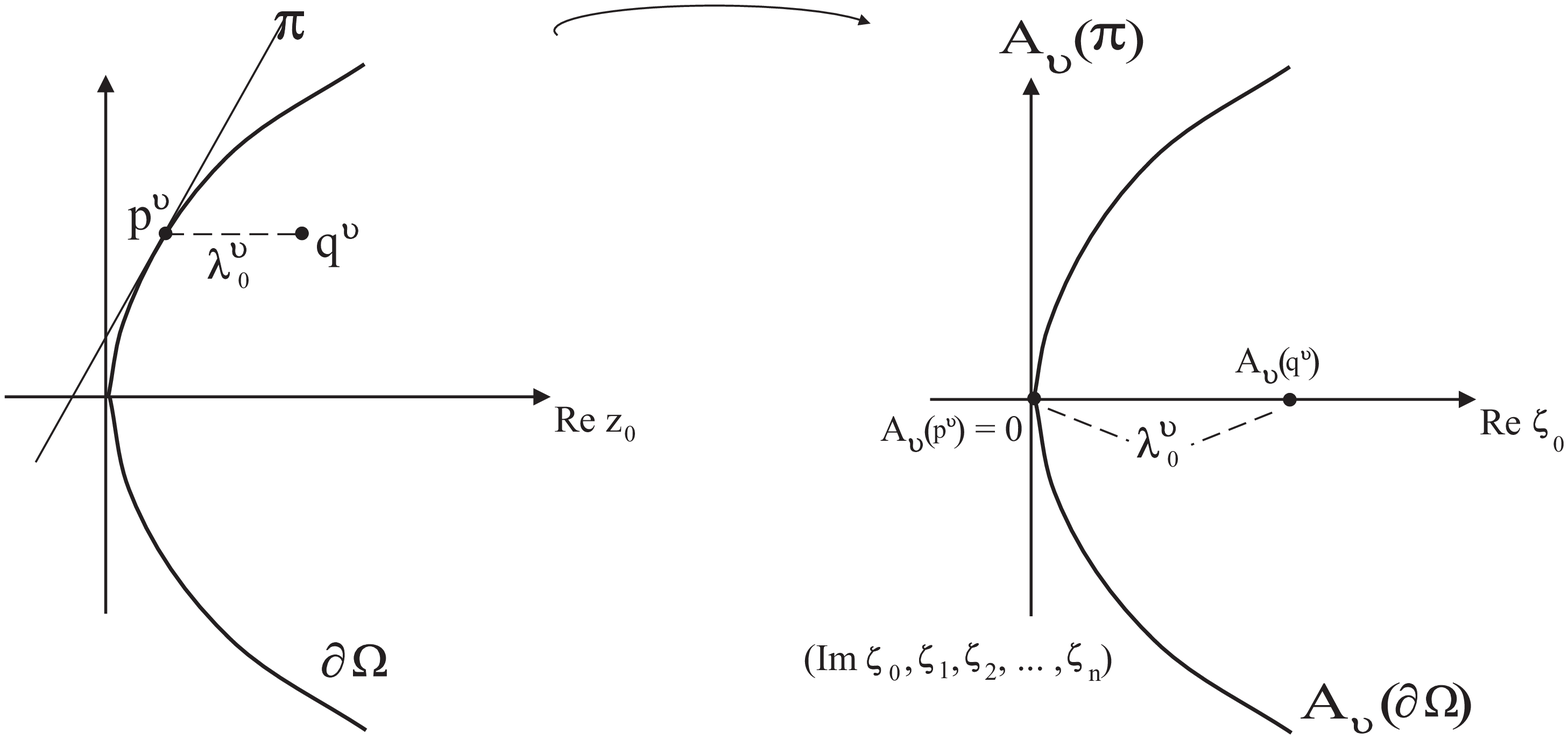}
\caption{The centering process.}
\end{figure}

Here the complex constants $\alpha_0^\nu, \ldots, \alpha_n^\nu$
are chosen so that $\alpha_0^\nu \to 1$ and $\alpha_m^\nu \to 0$
as $\nu \to \infty$ for $m=1,\ldots,n$, and such that the domain
$A_\nu (\Omega)$ is represented in a neighborhood of the origin by
a new $\mathcal C^k$ defining inequality
$$
\re \zeta_0 > \Psi_\nu (\im \zeta_0, \zeta_1, \ldots \zeta_2)
$$
satisfying
$$
\Psi_\nu (0, \ldots, 0) = 0 \hbox{ and } \nabla \Psi_\nu \big|_0 =
(0, \ldots 0).
$$

The next step is to consider a sequence of linear maps
$L_\nu:\CC^{n+1} \to \CC^{n+1}$ defined by
$$
L_\nu (\zeta_0, \ldots, \zeta_n) =
\left(\frac{\zeta_0}{\lambda_0^\nu},
\frac{\zeta_1}{\lambda_1^\nu}, \ldots,
\frac{\zeta_n}{\lambda_n^\nu}\right),
$$
where
$$
\lambda_0^\nu = q_0^\nu - p_0^\nu
$$
for each $\nu = 1, 2, \ldots$. The sequence of complex affine
mappings $\Lambda_\nu := L_\nu \circ A_\nu$ is the {\it dilatation
sequence} introduced\footnote{Pinchuk originally named it the {\it
stretching coordinates}.} by S. Pinchuk. The choice of
$\lambda_1^\nu, \ldots, \lambda_n^\nu$ is an important step in the
setup of the dilatation sequence. However, for the sake of smooth
exposition, it seems the best to postpone the explication until we
handle the Wong-Rosay Theorem (our first important application) a
bit later in the paper. It should be noted that Pinchuk's
dilatation sequence can be defined for a domain with boundary that
is not necessarily smooth near the reference boundary point $p$.

\subsection{Pinchuk's Scaling Sequence with Automorphisms}

Let $\Omega$ be a domain in $\CC^{n+1}$ with $\mathcal C^k$-smooth
boundary $\partial \Omega$, $k \ge 2$. Let $p$ be a boundary point
and $q^\nu$ be a sequence of points in $\Omega$ converging to $p$
as $\nu$ tends to infinity. Here we consider the important special
case when the sequence $q^\nu$ is given by $q^\nu = \varphi_\nu
(q)$, where $q$ is a point in $\Omega$ and where $\varphi_\nu$ is
a holomorphic automorphism of $\Omega$, for each $\nu=1,2,\ldots$.
Let the dilatation sequence $\Lambda_\nu: \Omega \to \CC^{n+1}$ be
as above, associated with the point sequence $q^\nu$. Then \it the
scaling sequence \rm introduced by S. Pinchuk is the following
sequence of maps:
$$
\sigma_\nu := \Lambda_\nu \circ \varphi_\nu : \Omega \to
\CC^{n+1}.
$$
Once the sequence $\varphi_\nu$ of automorphisms of $\Omega$ and a
point $q \in \Omega$ is given, the {\it orbit} $\varphi_\nu (q)$
and the affine adjustments $A_\nu$, which we call the {\it
centering maps}, are defined.  The only part of the scaling that
needs to be chosen is the sequence of dilating linear maps
$L_\nu$. The crux of the matter is to choose $L_\nu$ appropriately
so that
\begin{itemize}
\item[(\romannumeral1)] the $\sigma_\nu$ form a pre-compact normal
family, and
\item[(\romannumeral2)] a subsequential limit, say
$\widehat\sigma$, defines a holomorphic embedding of $\Omega$ into
$\CC^{n+1}$.
\end{itemize}
We shall see how such a simple idea produces significant results
in the subsequent sections.  On the other hand, it is not known
whether such a choice is always possible so that the scaling
sequence converges.

\subsection{Frankel's Scaling Sequence}

Before discussing the effect of the scaling method, it should be
mentioned that there is another way of constructing a scaling
sequence. With the same $\varphi_\nu$ and $q$ as above, S. Frankel
in his Ph.D. dissertation introduced the sequence
$$
\omega_\nu (z) \equiv [d\varphi_\nu (q)]^{-1} (\varphi_\nu (z) -
\varphi_\nu (q)).
$$
Notice that each $\omega_\nu$ embeds $\Omega$ into $\CC^{n+1}$.
>From the viewpoint of Pinchuk's scaling, one may see that the
differences between the two scaling methods are: {\bf (1)} that
the reference points $\varphi_\nu (q)$ are in the interior of
$\Omega$, and {\bf (2)} that the sequence $d\varphi_\nu (q)^{-1}$
replaces the role of centering followed by dilation. In the
scaling methods, the most delicate and important issues lie in the
convergence of the scaling sequence to a biholomorphic embedding
of the domain. One would also like to be able to guarantee that
the limit mapping is injective.  In the ensuing discussion we
shall invoke the notion of Kobayashi hyperbolicity. This is an
invariant version of the idea of boundedness.  A domain $\Omega$
or, more generally, a complex manifold is Kobayashi hyperbolic if
the Kobayashi distance (see \cite{KRA1}) is positive on $\Omega$.
The collection of Kobayashi hyperbolic manifolds is broad. For
example, any bounded domain is hyperbolic. Moreover, the upper
half-plane $\{z \in \CC \mid \Im z > 0\}$ and the Siegel upper
half-space $\{z \in \CC^{n+1} \mid \Im z_0 > |z_1|^2 + \ldots +
|z_n|^2 \}$ are Kobayashi hyperbolic.

\begin{theorem} {\rm (Frankel \cite{FRA})}  \sl
Let $\Omega$ be a convex, Kobayashi hyperbolic domain in
$\CC^{n+1}$. For any sequence $\varphi_\nu$ of automorphisms of
$\Omega$ and a point $q \in \Omega$, the sequence of maps defined
by
$$
\omega_\nu (z) = [d\varphi_\nu(q)]^{-1} (\varphi_\nu (z) -
\varphi_\nu (q))
$$
forms a pre-compact normal family, in the sense that its every
subsequence has a subsequence that converges uniformly on compact
subsets of $\Omega$.  Moreover, every subsequential limit is a
holomorphic embedding of $\Omega$ into $\CC^{n+1}$.
\end{theorem}

Notice that, from the viewpoint of this article at least, {\it
this theorem is most interesting when $\varphi_\nu (q)$
accumulates at a boundary point}. On the other hand, it turns out,
due to work of Kim and Krantz \cite{KIK1}, that Pinchuk's scaling
sequence can be selected to have the same conclusion in case the
domain is convex as in the hypothesis of this theorem.
Furthermore, the two scaling methods are indeed equivalent. More
precisely the following has been shown:

\begin{proposition} \sl
In addition to the hypothesis of Frankel's theorem above, assume
that the sequence $\varphi_\nu (q)$ accumulates at a boundary
point of $\Omega$ as $\nu \to \infty$. Let $\sigma_\nu$ denote
Pinchuk's scaling sequence.  Then we have the following
conclusion:
\begin{itemize}
\item[(\romannumeral1)] Every subsequence of $\sigma_\nu$ admits a
subsequence that converges uniformly on compacta to an injective
holomorphic mapping of $\Omega$ into $\CC^{n+1}$.
\item[(\romannumeral2)] Let $\widetilde \Omega$ denote the limit
domain of the Frankel scaling, and let $\widehat\Omega$ the limit
domain of the Pinchuk scaling.  Then these two domains are
biholomorphic to each other by a complex affine linear map.
\end{itemize}
\end{proposition}

\subsection{Normal Convergence of Sets}

We present the concept of Carath\'eodory kernel convergence of
domains which is relevant to the discussion of scaling methods and
general normal family of holomorphic mappings. For more detailed
discussions on this convergence, see p.76 of \cite{DUR}.

\begin{definition}[Caratheodory Kernel Convergence] \rm
Let $\Omega_\nu$ be a sequence of domains in $\CC^{n+1}$ such that
$\displaystyle{p \in \bigcap_{\nu=1}^\infty \Omega_\nu}$. If $p$
is an interior point of $\displaystyle{\bigcap_{\nu=1}^\infty
\Omega_\nu}$, the {\it Carath\'eodory kernel} $\widehat\Omega$
{\it at $p$} of the sequence $\{\Omega_\nu\}$ is defined to be the
largest domain containing $p$ having the property that each
compact subset of $\widehat\Omega$ lies in all but a finite number
of the domains $\Omega_\nu$. If $p$ is not an interior point of
$\displaystyle{\bigcap_{\nu=1}^\infty \Omega_\nu}$, the the
Carath\'eodory kernel $\widehat\Omega$ is $\{p\}$. The sequence
$\Omega_\nu$ of domains is said to {\it converge to its kernel at
$p$} if every subsequence of $\Omega_\nu$ has the same kernel at
$p$.
\end{definition}

We shall also say that a sequence $\Omega_\nu$ of domains in
$\CC^{n+1}$ {\it converges normally} if there exists a point
$\displaystyle{p \in \bigcap_{\nu+1}^\infty \Omega_\nu}$ such that
$\Omega_\nu$ converges to its Carath\'eodory kernel at $p$.

The motivation for this notion can be seen in the following
proposition.  We omit the proofs, as they are routine.

\begin{proposition}  \sl
Let $\Omega_\nu$ form a sequence of domains in $\CC^{n+1}$ that
converges normally to the domain $\widehat\Omega$.  Let $W \in
\CC^{m}$ be a Kobayashi hyperbolic domain. Then every sequence of
holomorphic mappings $f_\nu : \Omega_\nu \to W$ contains a
subsequence that converges uniformly on compacta to a holomorphic
mapping $\widehat f : \widehat\Omega \to \overline{W}$.
Furthermore, if $\{g_\nu : W \to \Omega_\nu\}$ forms a pre-compact
normal family, then every subsequential limit, say $\widehat g$,
has its image contained in the closure of $\widehat\Omega$.
\end{proposition}

It may be appropriate to remark that the topological set
convergence such as a version of local Hausdorff convergence can
replace the normal convergence in case the domains in
consideration are convex domains.

\section{Domains with Noncompact Automorphism Group}

\subsection{The theorem of Bun Wong and Jean-Pierre Rosay}


\begin{theorem} {\rm (Wong \cite{WON}, Rosay \cite{ROS})}  \sl
Let $\Omega$ be a bounded domain in $\CC^{n+1}$ with a sequence of
automorphisms $\varphi_\nu$ and a point $q \in \Omega$ such that
$\lim_{\nu \to \infty} \varphi_\nu (q) = p$ for some $p \in
\partial\Omega$.  If $\partial\Omega$ is $C^2$ strongly
pseudoconvex in a neighborhood $U$ of $p$, then $\Omega$ is
biholomorphic to the unit open ball $B$ in $\CC^{n+1}$.
\end{theorem}

\begin{theorem} {\rm (Wong \cite{WON})}  \sl
Every smoothly bounded domain in $\CC^{n+1}$ with transitive
automorphism group is biholomorphic to the unit open Euclidean
ball in $\CC^{n+1}$.
\end{theorem}
\medskip

\noindent \it Proof of the Wong-Rosay Theorem by the Method of
Scaling. \rm This proof is essentially due to S. Pinchuk.  It
consists of four typical steps for a scaling proof: {\bf (1)}
preparation, {\bf (2)} localization, {\bf (3)} dilatation, and
{\bf (4)} synthesis. \bigskip

\noindent \bf Step 1. Preparation. \rm Without loss of generality,
let us assume that $p$ is the origin $0$ in $\CC^{n+1}$. Since
$\Omega$ is strongly pseudoconvex at the origin, we may perform a
holomorphic coordinate change at the origin so that in an open
Euclidean ball $B(0, 10r)$ of radius $10r$ centered at the origin,
the set $\Omega \cap B(0,10r)$ can be defined by an
inequality\footnote{This is a concrete implementation of the
statement, discussed earlier, that a strongly pseudoconvex point
may be convexified by a biholomorphic mapping.}
$$
\rho(z) < 0
$$
where
$$
\rho(z) = - \re z_0 + |z_0|^2 + \ldots + |z_n|^2 + R(z)
$$
and (using Landau's notation)
$$
R(z)= o(|z_0|^2 + \ldots + |z_n|^2).
$$
In particular, choosing $r\ge0$ smaller if necessary, we may
arrange that
$$
|R(z)| \le \frac14 (|z_0|^2 + \ldots + |z_n|^2) ~~~ \forall z \in
B(0, 2r)
$$
and that the boundary $\partial\Omega$ is now strongly convex in
$B(0,2r)$.
\smallskip \\

It may be appropriate to remark at this juncture that a
$\calC^\infty$ smooth, strongly pseudoconvex boundary has a 4-th
order contact with a sphere. (See Fefferman \cite{FEF}, and also
\cite{HOR}.)
\medskip \\

\noindent \bf Step 2.  Localization. \rm As a consequence of the
preceding step, we see that there exists a holomorphic function
$h:B(0,2r) \to \CC$ such that
$$
h(0) = 1 \hbox{ and } |h(\zeta)|<1 \hbox{ for every } \zeta \in
\overline{\Omega} \cap B(0,2r) \setminus \{0\}.
$$
Now look at the automorphisms $\varphi_\nu$ of $\Omega$. Since
$\Omega$ is a bounded domain, every subsequence of $\varphi_\nu$
admits a subsequence that converges to a holomorphic mapping from
$\Omega$ into the closure $\overline{\Omega}$ of $\Omega$,
uniformly on compact subsets.  Let $\Phi$ be a subsequential limit
of a subsequence $\varphi_{\nu_k}$. Since $\Phi(q)=0$ and since
$\Phi:\Omega \to \overline{\Omega}$ is holomorphic, we may exploit
the uniform convergence on compact subsets to see that there
exists a relatively compact neighborhood $U$, say, of $q$ such
that $\varphi_{\nu_k} (U) \subset B(0,2r) \cap \Omega$ for all
sufficiently large $k$. Then consider $h \circ \varphi_{\nu_k} |_U
: U \to D$, where $D$ denotes the unit disc. This yields now that
$h \circ \Phi (0) = 1$. Hence the Maximum Modulus Principle
implies that $\Phi(z) = 0$ for every $z \in U$. Since $U$ contains
a non-empty open set, we conclude that $\Phi$ vanishes identically
in $\Omega$. We may now deduce, replacing $\{\varphi_\nu\}$ by one
of its subsequences, that for every compact subset $K$ of $\Omega$
there exists a positive integer $N$ such that
$$
\varphi_\nu (K) \subset B(0;r) \cap \Omega
$$
whenever $\nu \ge N$.
\medskip   \\

\noindent \bf Step 3. Dilatation. \rm Consider now the sequence
$\varphi_\nu (q)$ in $\Omega$. Let $q^\nu = \varphi_\nu (q)$ for
each $\nu$. Choose the boundary point $p^\nu$ as in the
construction of Pinchuk's dilatation sequence above.  Then choose
the centering map $A_\nu$ and the dilation map $L_\nu$ as above
for each $\nu$.  We may examine Pinchuk's dilatation sequence
$\Lambda_\nu := L_\nu \circ A_\nu$. It is a simple matter to check
that the sequence of domains $\Lambda_\nu (\Omega \cap U)$
converges normally to the domain
$$
V = \set{ (z_0, \ldots, z_n) \in \CC^{n+1} : \re z_0 \ge |z_1|^2 +
\ldots + |z_n|^2 }.
$$
Moreover, one may replace $\varphi_\nu$ by a subsequence again to
have that
$$
\Lambda_\nu (\Omega \cap U) \subset \mathcal{E} \qquad \forall \
\nu \ \hbox{ sufficiently large},
$$
where $\mathcal{E} = \set{ (z_0, \ldots, z_n) \in \CC^{n+1} : \re
z_0 \ge \frac12(|z_1|^2 + \ldots + |z_n|^2) }$.
\medskip   \\

\noindent \bf Step 4. Synthesis via Normal Families. \rm
Take a sequence $W_\nu$ of relatively compact subdomains of
$\Omega$ satisfying
$$
\overline{W_\nu} \subset W_{\nu+1} \hbox{ for every }
\nu=1,2,\ldots
$$
and
$$
\bigcup_{\nu=1}^\infty W_\nu = \Omega.
$$
Consider now the scaling sequence $\sigma_\nu = \Lambda_\nu \circ
\varphi_\nu$. Choosing a subsequence, we may assume that
$$
\varphi_\nu (W_\nu) \subset \Omega \cap B(0,r)
$$
for every $\nu=1,2,\ldots$. Then by the preceding step, the
scaling sequence $\sigma_\nu := \Lambda_\nu \circ \varphi_\nu
|_{W_\mu}$ forms a normal family for every $\mu$ as $\nu \to
\infty$. Notice also that, for every compact subset $K'$ of $V$,
the sequence $\sigma_\nu^{-1}$ maps $K'$ into $\Omega$ for
sufficiently large $\nu$. Altogether, one sees that any
subsequential limit of the scaling sequence becomes a
biholomorphic mapping from $\Omega$ onto $V$. Since $V$ is
biholomorphic to the unit open ball, the theorem is now proved.
\hfill  $\Box\;$
\bigskip

\subsection{Domains with Piecewise Levi-Flat Boundary that
Possess Noncompact Automorphism Group}

The main theorem of this section is the following:

\begin{theorem} {\rm (Kim-Krantz-Spiro \cite{KKS})}  \sl
Every generic analytic polyhedron in $\CC^2$ with noncompact
automorphism group is biholomorphic to the product of the unit
open disc in $\CC$ and a Kobayashi hyperbolic Riemann surface
embedded in $\CC^2$.
\label{Kim-Krantz-Spiro2005}
\end{theorem}

A clarification of some terminology is in order. By an {\it
analytic polyhedron} in $\CC^{n+1}$, we mean a bounded domain, say
$\Omega$ in $\CC^{n+1}$, admitting an open neighborhood $U$ of the
closure $\overline{\Omega}$ of $\Omega$ and a finite collection of
holomorphic functions $f_j:U \to \CC$, $j=1,\ldots,N$, such that
$$
\Omega = \set{ z\in \CC^{n+1} : |f_1(z)|<1, \ldots, |f_N(z)|<1 }.
$$
The collection $\set{f_1, \ldots, f_N}$ is usually called a
defining system for $\Omega$. The choice for defining system is
not unique. An analytic polyhedron is called {\it generic} (or,
{\it normal}), if it admits a defining system $\set{f_1, \ldots,
f_N}$ satisfying the following additional condition:
$$
df_{i_1} |_p \wedge \cdots \wedge df_{i_k} |_p \neq 0
$$
whenever the condition $|f_{i_1}(p)| = \ldots = |f_{i_k}(p)| = 1$
holds for any un-repeated indices $i_1, \ldots, i_k \in \set{1,
\ldots, N}$. We restrict our attention as usual to the {\it
bounded} analytic polyhedra. By a theorem of H. Cartan mentioned
earlier, the automorphism group of our analytic polyhedron is a
finite-dimensional Lie group. The non-compactness of the
automorphism group is therefore equivalent to the existence of a
sequence $\varphi_\nu \in \aut\Omega$ and a point $q \in \Omega$
such that the point sequence $\varphi_\nu (q)$ accumulates at a
boundary point.
\bigskip

Notice that Theorem \ref{Kim-Krantz-Spiro2005} improves the
following results:

\begin{theorem} {\rm (Kim-Pagano \cite{KIP})}
\label{Kim-Pagano2000} \sl
Let $\Omega$ be a generic analytic polyhedron in $\CC^2$ with
noncompact automorphism group.  Then the holomorphic universal
covering space of $\Omega$ is biholomorphic to the bidisc.
\end{theorem}

It is worth mentioning the following theorem. It concerns only
convex analytic polyhedra, but is valid in all dimensions.

\begin{theorem} {\rm (Kim \cite{KI1})} \label{Kim1992}  \sl
Every convex generic analytic polyhedron in $\CC^{n+1}$ with
noncompact automorphism group is biholomorphic to the product of
the unit open disc in $\CC$ and a convex domain in $\CC^{n}$. In
particular, in case $n=1$, the product domain is biholomorphic to
the bidisc.
\end{theorem}
\bigskip

\noindent \bf Sketch of proof. \rm
We treat all three theorems simultaneously. Let $\Omega$ be a
generic analytic polyhedron in $\CC^{n+1}$ with noncompact
automorphism group. Then there exist a boundary point $p \in
\partial\Omega$, an interior point $q \in \Omega$ and a sequence
$\varphi_j \in \aut\Omega$ such that $\lim_{j\to\infty} \varphi_j
(q) = p$. Since the boundary of $\Omega$ is piecewise Levi flat,
we divide the proof of each theorem above into the following two
cases: {\bf (1)} the case that boundary $\partial\Omega$ is
singular at $p$, and {\bf (2)} the case that boundary
$\partial\Omega$ is smooth and Levi flat in a neighborhood of $p$.
\medskip   \\

\noindent \underbar{\it Case 1.} The boundary $\partial\Omega$ is
singular at $p$.

\medskip

Recall that our domain is a generic analytic polyhedron. In
complex dimension 2, therefore, it is possible to choose a
defining system $f_1, \ldots, f_N$ such that there exist exactly
two functions, say $f_1, f_2$ (shuffling the indices if
necessary), such that $|f_1(p)|=|f_1(p)|=1$ with $df_1|_p \wedge
df_2|_p \neq 0$. Hence it is simple to realize that there exists a
plurisubharmonic\footnote{A real-valued continuous function
$\psi:\Omega \to \RR$ defined in a domain $\Omega$ in $\CC^n$ is
called {\it plurisubharmonic} if it is subharmonic when restricted
to any complex affine line. See \cite{KRA1}.}
function $\psi$ defined in an open neighborhood of the closure of
$\Omega$ such that
$$
\psi (p)=0, \hbox{ and } \psi (z) < 0 \hbox{ for every } z \in
\overline{\Omega} \setminus \{p\}.
$$
Such a function is called a {\it plurisubharmonic} (psh for short)
{\it peak function} for $\Omega$ at $p$. The maximum modulus
principle immediately implies in particular that there are no
non-trivial analytic varieties in $\partial\Omega$ passing through
$p$.  Moreover, it is known that the following localization
principle holds (see \cite{BER}, and also \cite{BYGK} for detailed
arguments, for instance):

\medskip

\it Let $U$ be an open neighborhood of $p$. For every compact
subset $K$ of $\Omega$, there exists a positive integer $j_K$ such
that $\varphi_j (K) \subset U$ for every $j \ge j_K$. \rm
\medskip

Now, using the mapping $(f_1, f_2):U \to \CC^2$ followed by a
linear fractional mapping of $\CC^2$, one can construct a
biholomorphism-into $\Psi: U \to \CC^2$ with $\Psi(p)=(0,0)$ and
$\Psi(U \cap \Omega) = \Psi(U) \cap V^2$ where
$$
V^2 = \set{(z_0, z_1) \in \CC^2 : \re z_0 \ge0, \re z_1\ge 0}.
$$
\medskip
Let us write $\Psi (\varphi_j (q)) = (t_{j,0}, t_{j,1})$. Then
consider the dilatation sequence
$$
\Lambda_j (z_0, z_1) = \left( \frac{z_0 - \im t_{j,0}}{\re
t_{j,0}}, \frac{z_1 - \im t_{j,1}}{\re t_{j,1}} \right).
$$
Finally, consider the scaling sequence
$$
\sigma_j := \Lambda_j \circ \Psi \circ \varphi_j
$$
for $j = 1,2,\ldots$. It still requires some checking, but it
follows that a subsequential limit of this sequence gives rise to
a biholomorphic mapping from $\Omega$ onto the domain $H^2$ which
is in turn biholomorphic to the bidisc. Thus Theorems 4.3 and 4.4
are proved in this case. For Theorem 4.5, one takes into
consideration that $\Omega$ is a convex domain in $\CC^{n+1}$.
Note that every complex analytic variety contained in a convex
Levi-flat hypersurface is an open subdomain of a complex affine
hyperplane. Now choose a defining system $f_0, \ldots, f_N$ so
that
$$
f_0(p) = \ldots = f_k (p)=1, |f_{k+1}(p)|<1, \ldots, |f_N (p)| < 1
$$
and
$$
df_1|_p \wedge \cdots \wedge df_k|_p \neq 0.
$$
The maximal variety in the boundary of $\Omega$ passing through
$p$ is represented by
$$
V_p = \set{z \in \CC^{n+1} : f_0(z) = \ldots = f_k (z)=1,
|f_{k+1}(z)|<1, \ldots, |f_N (z)| < 1}.
$$

Notice that $\dim_\CC V_p = n-k$. It is possible that $n=k$, and
consequently that $V_p$ is a single point. But it is always the
case that $k\ge 0$. Then one may use the mappings $f_0, \ldots,
f_k$ to see that there exists an open neighborhood, say $W$, of
$V_p$ such that there exists a holomorphic embedding $\Psi : W \to
\CC^{n+1}$ such that
\bigskip

\begin{itemize}
\item[(a)] $\Psi (V_p) = \set{(0, \ldots, 0; z_{k+1}, \ldots, z_n)
\in \CC^{n+1} : (z_{k+1}, \ldots, z_n) \in \Omega'}$, where
$\Omega'$ is a convex domain containing the origin in $\CC^{n-k}$,
and
\item[(b)] $\Psi (W \cap \Omega) = \Psi (W) \cap U$ where
$$
U = \set{z \in \CC^{n+1} : \re z_0 \ge 0, \ldots, \re z_k \ge 0,
(z_{k+1}, \ldots, z_n) \in \Omega'}.
$$
\end{itemize}

\bigskip
\noindent Write $\Psi \circ \varphi_j (q) = (t_{j,0}, \ldots,
t_{j, n})$ and then consider the dilatation map
$$
\lambda_j (z) = \left( \frac{z_0 - \im t_{j,0}}{\re t_{j,0}},
\ldots, \frac{z_k - \im t_{j,k}}{\re t_{j,k}}; z_{k+1}, \ldots,
z_n \right).
$$
Again it follows that a subsequential limit of the sequence
$$
\sigma_j := \Lambda_j \circ \Psi \circ \varphi_j
$$
gives rise to a biholomorphic mapping from $\Omega$ onto $U$ which
in turn is biholomorphic to the product of a $k+1$ dimensional
polydisc and the domain $\Omega'$ in $\CC^{n-k}$. This proves
Theorem 4.5 in the present case.
\bigskip  \\

\underbar{\it Case 2.} The boundary $\partial\Omega$ is smooth and
Levi flat at $p$. This case is easier to handle when $\Omega$ is
convex, even when $\dim_{\CC^{n+1}} \Omega = n+1$ for an arbitrary
non-negative integer $n$. Consider the maximal variety $V_p$ in
the boundary $\partial\Omega$ passing through $p$. As argued
earlier, it follows that $V_p$ is a convex, open subset of a
complex affine hyperplane. By a complex affine linear change of
coordinates, say by an affine linear biholomorphism
$\psi:\CC^{n+1} \to \CC^{n+1}$, we may assume that $\psi(p)=0$,
that $\psi(V_p) \subset \set{(z_0, \ldots, z_n) : z_0 = 0}$ and
that the domain $\Omega$ at $p$ is contained in the half-space
defined by the inequality $\re z_0 \ge 0$. Now we apply the
scaling method as before. Let $\psi\circ\varphi_j (q) = (t_{j,0},
\ldots, t_{j,n})$.  Then define the dilatation mapping by
$$
\Lambda_j (z_0, \ldots, z_n) = \left( \frac{z_0 - \im t_{j,0}}{\re
t_{j,0}}, z_1, \ldots, z_n \right).
$$
Then it turns out that the scaling sequence
$$
\sigma_j : = \Lambda_j \circ \psi \circ \varphi_j
$$
yields a subsequential limit which becomes a biholomorphic mapping
from $\Omega$ onto the product of the upper half plane in $\CC$
and the $n$-dimensional convex domain $\CC^{n+1}$. This completes
our sketch of the proof to Theorem 4.5.
\medskip    \\

In case the analytic polyhedron is merely generic, and not
necessarily convex, the situation is much more complicated. Thus
it is natural that one focuses on the case of complex dimension
two.

Then the maximal variety $V_p$ is a Riemann surface that is
Kobayashi hyperbolic. The uniformization theorem of Riemann
surface theory yields a holomorphic covering map $\pi:D \to V_p$
from the open unit disc $D$ in $\CC$ onto $V_p$.  Then extend it
trivially to the map $\widetilde\pi (z_1, z_2) = (\pi(z_1), z_2)$.
Since the normal bundle for $V_p$ in $\CC^2$ is trivial, one may
take an open neighborhood $U$ for $V_p$ in $\CC^{n+1}$ in such a
way that $\Omega \cap U$ is connected and that $\widetilde\pi$
gives rise to a local biholomorphism, say $\widehat\pi$, from an
open neighborhood of $D \times \set{0}$ onto $U$.  One can arrange
also that $\widetilde\pi(0,0) = p$.

Then consider the domain ${\widetilde\Omega}_{\rm loc}$ which is a
connected component of $\widehat\pi^{-1} (\Omega \cap U)$
containing the origin.  Then take a lifting of the sequence
$\varphi_j(q)$ via $\widetilde\pi$. There are many liftings.
Choose therefore one that converges to the origin $(0,0)$ for
instance. Now build a dilatation mapping $\Lambda_j$ for
$\widetilde\Omega$, formally the same as in the convex case (with
an adjustment; see \cite{KKS} for details), with respect to the
sequence chosen here. Then it turns out that the sequence of
mappings
$$
\varphi_j^{-1} \circ \widehat\pi \circ \Lambda_j
$$
yields a subsequential limit, say $\Psi$, from the product $D
\times H$ of the open unit disc $D$ and the half-plane $H = \set{z
\in \CC : \re z \ge 0}$ onto our generic analytic polyhedron
$\Omega$. Using normal families arguments, it is not hard to
deduce that $\Psi$ is a holomorphic mapping with its Jacobian
vanishing nowhere on $D \times H$.  Moreover, it turns out that
this map preserves the Kobayashi-Royden infinitesimal metric.
Therefore it preserves the Wu metric (see below) as well.
\medskip

The Wu metric (see \cite{WU}) here can be quickly understood as
follows. At each point $p$ of a Kobayashi hyperbolic domain $G$ in
$\CC^{n+1}$, consider the tangent space $T_p G$ and the collection
of vectors with Kobayashi length not exceeding 1. This set is
sometimes called the Kobayashi indicatrix at $p$. Endow $T_p G$
with an arbitrarily chosen Hermitian inner product. (This choice
of Hermitian inner product is neither unique nor natural, but it
will not cause any problems at the end.) Then consider the
ellipsoids, say $E=E_H$, in $T_p G (= \CC^{n+1})$ given by
$$
E_H = \{ v \in \CC^{n+1} : v^* H v \le 1 \}
$$
where $H$ is a positive definite Hermitian $(n+1)\times(n+1)$
matrix and where $v^*$ denotes the conjugate transpose of $v$.
Denote by $\mathcal{Q}$ the set of such $E_H$ containing the
Kobayashi indicatrix. Then take $1/(\det H)$ as the volume of
$E_H$. Then it turns out that the element in $\mathcal{Q}$ with
the smallest volume is uniquely determined, regardless of the
choice of the Hermitian inner product on $T_p$. See \cite{WU}. Now
let this minimum volume ellipsoid define a Hermitian inner
product, say $h_p$, on $T_p G$. The assignment $p \mapsto h_p$
defines the Wu metric on $G$.  It is immediate from the invariance
of the Kobayashi metric that the Wu metric is invariant under
biholomorphic maps. It is obviously Hermitian, in the sense that
it defines a Hermitian inner product on each tangent space.  It
has been shown that $p \mapsto h_p$ is $\mathcal{C}^0$
(continuous) in general.
\medskip

Now since the Wu metric $h_{D\times D}$, say, of the bidisc $D
\times D$ is real-analytic, so are the Wu metrics $h_{D \times H}$
of $D \times H$ and $h_\Omega$ of $\Omega$, respectively, since
$\Psi_* h_{D \times H} = h_\Omega$. Since the Kobayashi metric of
$\Omega$ is complete, so is the Wu metric $h_\Omega$. At this
point, one may apply the proof of the Cartan-Hadamard Theorem in
Riemannian Geometry to conclude that $\Psi: D \times H \to \Omega$
is indeed a covering mapping.  This yields Theorem
\ref{Kim-Pagano2000}.
\medskip  \\

Finally, for Theorem \ref{Kim-Krantz-Spiro2005}, one has to
analyze the covering mapping $\Psi$ as well as its deck
transformation group more precisely. Although we are omitting the
details here, it should not be difficult for the reader to see
that the constructions for $\widetilde\pi$ as well as
$\widehat\pi$ will reflect the nature of the covering map $\pi:D
\to V_p$ without any essential changes. Hence it was shown in
\cite{KKS} through a careful analysis that indeed $\Omega$ is
biholomorphic to the product of the open unit disc and the maximal
variety $V_p$, and that the deck transformation group
$\Gamma_\Psi$ for the covering mapping $\Psi: D \times H \to
\Omega$ is in fact $\Gamma_\pi \times \set{\rm id}$, where
$\Gamma_\pi$ denotes the deck-transformation group for the
uniformization map $\pi:D \to V_p$.  This is how Theorem
\ref{Kim-Krantz-Spiro2005} follows. \hfill  $\Box\;$
\bigskip  \\

Notice that this analysis gives a rather complete classification
of complex two-dimensional generic analytic polyhedra that possess
noncompact automorphism group. Hence it seems appropriate to pose
the following question here.

\begin{Problem} \sl Classify all non-generic analytic polyhedra
with noncompact automorphism group.
\end{Problem}

The main difficulty in this problem seems to lie in the question
of how to adjust the scaling at a singular boundary point.

\subsection{A Digression on Finite Type}

It is a straightforward calculation (see \cite{KRA1}) to see that
a strongly pseudoconvex boundary point $P$ is flat to order 2.
That is to say, the maximum possible order of contact of the
boundary at $P$ with a one-dimensional complex analytic variety is
2. See Figure 4. It is useful in this subject to be able to
generalize this concept.  For simplicity in this subsection we
restrict our attention to two-dimensional complex space.  The
entire story in all dimensions is sketched out in \cite{KRA1}.

\begin{figure}
\centering
\includegraphics[height=2.65in, width=2.75in]{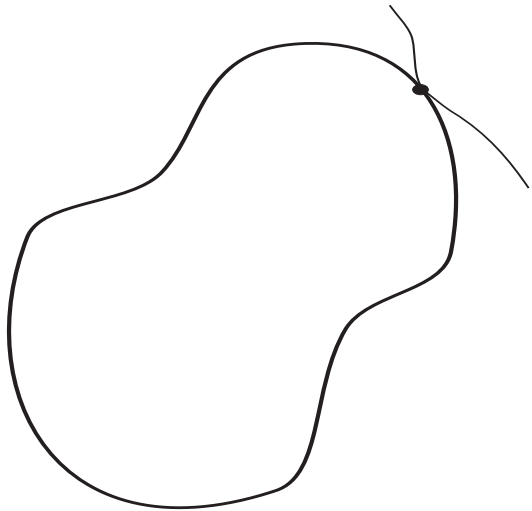}
\caption{Order of contact at a strongly pseudoconvex point.}
\end{figure}

Let $\Omega = \{z \in \CC^{n+1}: \rho(z) < 0\}$ be a smoothly bounded
domain and fix a point $p \in \partial \Omega$.  Let $V$ be a
nonsingular, one-dimensional complex analytic variety that passes
through $p$.  Then the {\it order of contact} of $V$ with $\bO$ at
$p$ is the greatest positive integer $k$ such that
$$
|\rho(z)| \leq C \cdot |z - p|^k
$$
for $z \in V$ near $p$ and some constant $C \ge 0$. We say that
$P$ is of {\it finite geometric type} in the sense of
Catlin/D'Angelo/Kohn if there is an upper bound $m$ on the order
of contact of analytic varieties with $\bO$ at $p$.  The least
such integer $m$ is called the {\it type} of the point $p$.  As
previously noted, a strongly pseudoconvex point is of type 2 (see
\cite{KRA1} for the details). As a very simple illustrative
example, for $k$ a positive integer let
$$
E_k = \{z = (z_1, z_2) \in \CC^2: |z_1|^2 + |z_2|^{2k} < 1\} \, .
$$
Then one may calculate (again see \cite{KRA1}) that any boundary
point of the form $(e^{i\theta},0)$ is of finite type $2k$. The
notion of ``type'' is a means of measuring the flatness of the
boundary in a complex analytic sense. One of the most important
facts about finite type in complex dimension 2 is that the
geometric definition given here is equivalent to an analytic
definition involving commutators of vector fields.  To wit, we may
assume by a normalization of coordinates that $\partial
\rho/\partial z_1(p) \ne 0$.  Define the vector field
$$
L = \frac{\partial \rho}{\partial z_1}(p) \frac{\partial}{\partial
z_2} - \frac{\partial \rho}{\partial z_2}(p)
\frac{\partial}{\partial z_1} \, .
$$
Then $L$ is a tangential holomorphic vector field near $p$
(because $L \rho = 0$). A {\it first-order commutator} is, for us,
an expression of the form $[L, \overline{L}] = L \overline{L} -
\overline{L} L$. A {\it second-order commutator} is the commutator
of $L$ or $\overline{L}$ with a first-order commutator. And so
forth. We say that $p$ is of {\it analytic type} $m$ if any
commutator of order not exceeding $m - 1$ has no complex normal
component but that some commutator of order $m$ {\it does} has a
complex normal component.  It is a fundamental result of Kohn
\cite{KOH} (generalized to higher dimensions by Bloom and Graham
\cite{BLG}) that the boundary point $p$ is of finite geometric
type $m$ if and only if it is of finite analytic type $m$.

One immediate consequence of this
characterization is the semicontinuity of type:  If $p \in
\partial \Omega$ is a point of finite type $m$ then there is a
small boundary neighborhood $U$ of $p$ so that all point of $U$
are of finite type not exceeding $m$.  In complex dimensions
greater than 2, this semicontinuity (as stated here) fails.  But
there is a substitute result that was proved by D'Angelo in
\cite{DAN1}.  In any event, it follows from these results that if
$\Omega$ is smoothly bounded and if each point of $\bO$ is of
finite type, then there is an upper bound $M$ so that the type of
{\it every} boundary point does not exceed $M$. In deep work
\cite{DIF}, Diederich and Forn\ae ss showed that any domain with
real analytic boundary, in any complex dimension, is of finite
type. Thus domains with real analytic boundary form an important
class of examples in this subject. (We recommend the readers to
read an alternative proof by D'Angelo in \cite{DAN2}.) It is
perhaps instructive to contrast such a domain with a boundary that
is Levi flat.  Such a boundary is foliated by complex analytic
varieties, so that one sees immediately that each boundary point
is of infinite type. The provenance of the concept of finite type
was the study of the $\dbar$-Neumann problem (see \cite{KOH}).
Since that time, finite type has assumed a rather prominent
position in function theory, mapping theory, and related areas.
See \cite{DAN1}, \cite{DAN2} for a full account of this central
idea.

\subsection{A Theorem of Bedford-Pinchuk}

In the preceding sections, we discussed domains with noncompact
automorphism group in the extreme cases when the boundaries are
either strongly pseudoconvex or Levi flat.  The intermediate
concept, called the boundary of finite type in the sense of
Catlin/Kohn/D'Angelo, encompasses a large class of weakly
pseudoconvex domains with smooth boundary.  We now present the
following theorem pertaining to this class.

\begin{theorem} {\rm (Bedford-Pinchuk \cite{BEP1})}    \sl
A bounded pseudoconvex domain $\Omega \ss \CC^2$ having real
analytic boundary and admitting noncompact automorphism group is
biholomorphic to
$$
E_m = \{(z,w) \in \CC^2 : |z|^2 + |w|^{2m} < 1 \}
$$
for some positive integer $m$.
\end{theorem}

\noindent In fact some unpublished remarks of David Catlin show
that the hypothesis ``real analytic'' may be weakened to ``finite
type''. It should be also mentioned that Bedford and Pinchuk have
extended the above theorem to broader classes of domains. Before
discussing the proof of this theorem, we remark that several
generalizations of this theorem have been carried out, mostly by
Bedford and Pinchuk. See \cite{BEP1}, \cite{BEP2}, \cite{BER} for
precise results.  In particular Bedford and Pinchuk can prove
results in higher dimensions with certain restrictions such as
convexity. However, we shall focus on the original proof, as it
reflects the essential methods for this case.
\medskip \\

\noindent \bf Sketch of the Proof of Theorem 4.7. \rm
The noncompactness of $\aut\Omega$ again implies the existence of
$\varphi_j \in \aut\Omega$, $p \in \partial\Omega$ and $q \in
\Omega$ such that $\lim_{j\to\infty} \varphi_j (q) = p$. Then one
may choose a holomorphic local coordinate system $(z,w)$ such that
$p$ becomes the origin, and such that there exists an open
neighborhood $U$ of the origin in which the domain $\Omega$ is
represented by the inequality
$$
\re w \ge H(z) + R(z,w)
$$
where:
\smallskip \\
\noindent {\bf (i)} $H(z)$ is a homogeneous subharmonic polynomial
in $z, \bar z$ of degree $2m$ without harmonic terms. Here $m$ is
a positive integer, and
\smallskip \\
\noindent {\bf (ii)} $R(z,w) = o(|z|^{2m} + |\im w|)$.
\medskip  \\
Then it follows by a careful application of the scaling method
that $\Omega$ is biholomorphic to the domain
$$
M(\Omega, p) = \set{(z,w) \in \CC^2 : \re w \ge H(z)} \, .
$$
For a detailed argument a useful reference other than the original
paper by Bedford-Pinchuk is the theorem on p.\ 620 of \cite{BER}
by Berteloot. At this juncture, we simply use the fact that the
mapping $(z,w) \mapsto (z, w+it)$ is an automorphism of
$M(\Omega,p)$ for every $t \in \RR$. This produces a noncompact
one-parameter subgroup of automorphisms, say $\psi_t$, for
$\Omega$. But then Bedford and Pinchuk showed that there exists a
boundary point $p' \in \partial\Omega$ such that
$$
\lim_{t \to \pm \infty} \psi_t (z,w) = p'
$$
for every $(z,w) \in \Omega$.  Furthermore they show that this
turns into a smooth parabolic holomorphic vector field action at
$p'$ on the boundary of $\Omega$.
The next step is to establish that this parabolic orbit by this
one-parameter family of automorphisms is not too tangent to the
boundary; in fact it is shown in [BEP1] that the
order of contact is at most quadratic.
Then we use this parabolic orbit to scale $\Omega$ again. This
will again yield that the original domain $\Omega$ is
biholomorphic to the domain $M(\Omega,p')$ defined by
$$
\re w \ge \Psi(z) \, ,
$$
where $\Psi$ is a homogeneous subharmonic polynomial in $z, \bar
z$ of an even degree, say $2\ell$.  But this time the important
thing is that the one-parameter action from $\Omega$ survives in
the domain $M(\Omega,p')$.  This introduces a further restriction
to $\Psi(z)$ and eventually forces that $\Psi (z) = |z|^{2\ell}$,
yielding the desired conclusion.
\endpf
\bigskip  \\
As one can see from the proof, the assumption that $\Omega$ has
global smoothness (indeed, real analyticity) and finite type is
essential, as one does not know where $p'$ will be located in
$\partial\Omega$. Attempts to obtain the same conclusion from the
weaker assumption that $\partial\Omega$ is real analytic of finite
type at the initial orbit accumulation point $p$ cannot be
successful, as there are obvious counterexamples such as the one
defined by $\re z_0 > |z_1|^8 + 2 \re \bar z_1^3 z_1^5
$, for
instance.

\section{The Greene-Krantz Conjecture}

The classification program for domains with noncompact
automorphism group is far from being complete, even for the case
of very smooth, or piecewise smooth, boundaries. On the other hand,
it seems natural at this juncture to mention the following
outstanding conjecture by Greene and Krantz:

\begin{conjecture} \sl
Let $\Omega$ be a bounded domain in $\CC^{n+1}$ with  $\mathcal
C^\infty$ boundary. If there exists a sequence $\varphi_\nu$ of
automorphisms of $\Omega$ and a point $q \in \Omega$ such that the
orbit $\varphi_\nu (q)$ accumulates at a boundary point $p$ of
$\Omega$, then $p$ is of finite type in the sense of
D'Angelo/Catlin/Kohn.
\end{conjecture}

The full conjecture is still open. The purpose of this section is
to introduce some partial results supporting the conjecture,
discovered by means of the scaling method.
The first partial result we mention here is the following
reformulation of Theorem \ref{Kim1992}:

\begin{proposition} \sl
Let $\Omega$ be a bounded, convex domain in $\CC^{n+1}$ with a
boundary point $p \in \partial\Omega$ admitting an open
neighborhood $U$ such that $\partial\Omega \cap U$ is Levi flat at
every point. Then no automorphism orbit of $\Omega$ can accumulate
at $p$.
\end{proposition}

The basis for this proposition is as follows: If there were an
automorphism orbit accumulating at $p$, then Theorem\ref{Kim1992}
from the hypothesis implies that $\Omega$ is biholomorphic to the
product of a convex domain and the open unit disc. Therefore,
$\Omega$ is a (trivial) fiber space over $\Omega$. A theorem of A.
Huckleberry (\cite{HUC}) says that this cannot be biholomorphic to
a bounded domain with a strongly pseudoconvex boundary point.
Since any bounded domain with entirely smooth boundary must
admit\footnote{Consider the function $f(x) = \|x\|$ that represent
the Euclidean distance between the origin and the point $x \in
\overline{\Omega}$. Since $\overline{\Omega}$ is compact, the
function $f(x)$ assumes the maximum, at a boundary point $p$, say,
of $\Omega$. Then $\partial\Omega$ has a sphere contact at $p$ so
that the whole domain $\Omega$ is included in the sphere. Then $p$
is in fact a strongly convex (hence, strongly pseudoconvex)
boundary point.} a strongly pseudoconvex boundary point, it leads
us to a contradiction. Thus the proposition follows immediately.
\hfill $\Box$
\bigskip

It should be observed that an infinite type boundary point need
not admit a neighborhood in which every boundary point is of
infinite type (consequently Levi-flat). A primary example is given
by the origin for the domain defined by
$$
\re w \ge \exp \left( -\frac1{|z|^2} \right).
$$
Indeed, Greene and Krantz demonstrated the following, when they
posed the aforementioned conjecture:
\begin{proposition} \sl
The automorphism group of the domain in $\CC^2$ defined by
$$
|z|^2 + 2 \exp \left( - |w|^{-2} \right) < 1
$$
is compact. In particular, there is no automorphism orbit
accumulating at any boundary point $(e^{i\theta}, 0)$, of infinite
type.
\end{proposition}

The original proof of this proposition exploited the fact that the
domain in consideration is Reinhardt, admitting full rotational
symmetry. However, it turns out that, for the purpose pertaining
to the Greene-Krantz conjecture, the obstruction against the
existence of automorphism orbits accumulating at the point of such
exponential infinite type boundary point is purely local.
Consider the following result.

\begin{theorem} \rm (Kim/Krantz \cite{KIK1}) \sl
Let $\Omega$ be a domain in $\CC^2$ with a boundary point $p$
which admits an open neighborhood $U$ and an injective holomorphic
mapping $\Psi:U \to \CC^2$ such that $\Psi(p)=(0,0)$ and
$$
\Psi (U \cap \Omega) = \set{(z,w) \in \Psi(U) : \re w \ge
\psi(|z|)}
$$
where $\psi:\RR \to \RR$ is a $C^\infty$ smooth function
satisfying:
\begin{itemize}
\item[(\romannumeral1)] $\psi$ is $C^\infty$ smooth.
\item[(\romannumeral2)] $\psi (t) = 0, ~\forall t\le 0$, and
$\psi''(t) \ge 0, ~\forall t\ge0$.
\item[(\romannumeral3)] $\psi (t) = \exp (- \mu(t)^{-1})$ for some
$\mu(t)$ that is a non-negative smooth function vanishing to a
finite order at $t=0$.
\end{itemize}
Then there is no holomorphic automorphism orbit of $\Omega$
accumulating at $p$.
\end{theorem}

\noindent \bf Sketch of the Proof. %
\rm The detailed argument of the proof given in \cite{KIK1} is
long and tedious. On the other hand, the key ideas are as follows.
Expecting a contradiction, assume to the contrary that there is an
automorphism orbit $\varphi_\nu (q)$ converging to $p$. Then apply
the scaling technique to the domain $\Omega$. Calculations show
that the scaled limit domain is biholomorphic to one of the
following domains:
\begin{itemize}
\item[(\romannumeral1)] the open unit ball $\BB$ in $\CC^2$,
\item[(\romannumeral1)] the open unit bidisc $D = \{(z,w) \in
\CC^2 : |z|<1, |w|<1 \}$
\item[(\romannumeral1)] the domain $T = \{(z,w) \in \CC^2 : \Re
z \ge \exp (\Re w) \}$.
\end{itemize}
These three domains occur depending upon the tangency of the orbit
$\varphi_\nu (q)$ to the boundary. If the orbit is very tangential
to the strongly pseudoconvex part of the boundary, then the first
possibility appears. If the orbit is not so tangential to the
boundary, than the bidisc shows up as the limit domain of the
scaling process. The appropriate intermediate exponential tangency
of the orbit to the strongly pseudoconvex part of the boundary
produces the 3rd possibility. Now notice that the convergence of
scaling (see \cite{KIK1}) implies that the original domain
$\Omega$ has to be biholomorphic to one of the domains listed
above. But none of these possibilities can occur.  In the first
case, the domain $\Omega$ should be homogeneous, as the ball is.
Then one may choose a non-tangential sequence of automorphism
orbit accumulating at $p$. Then scaling will show that the domain
$\Omega$ is biholomorphic to the bidisc.  This shows that the the
domain $\Omega$ is biholomorphic to both the ball and the bidisc.
This contradicts the theorem of Poincar\'e which says that there
does not exist any biholomorphism between the ball and the bidisc
in complex dimension two.  For the second case, a mirror-image
argument implies the same kind of contradiction. The third case is
much the same. Since the domain $T$ contains a real 3-dimensional
subgroup without fixed points, one finds an automorphism orbit of
$\Omega$ that is either non-tangential to the boundary or very
tangential to the line of boundary points of infinite type.  In
either case, one gets the bidisc as the new scaled limit.  Then
one arrives at a contradiction as before. This, altogether,
contradicts the theorem of convergence of Pinchuk's scaling method
for the convex case (See \cite{KIK1}, stated as Proposition 3.2 in
Section 3.3.). Therefore one is led to the conclusion that there
are no automorphism orbits in $\Omega$ accumulating at $p$, as
claimed. \hfill $\Box$
\bigskip  \\

Digressing slightly, it is worth observing the following results
of J. Byun (\cite{BY1}, \cite{BY2}):

\begin{theorem} {\rm (Byun)} \sl Let $\Omega$ be a domain in
$\CC^2$. Assume that there exists a point $p \in \partial\Omega$
admitting an open neighborhood $U$ in $\CC^2$ satisfying the
conditions
\begin{itemize}
\item[(1)] the boundary $\partial\Omega$ is $\calC^\infty$,
pseudoconvex, and of finite type in $U$, and
\item[(2)] the D'Angelo finite type of $\partial\Omega$ at $p$ is
strictly greater than that of other points in $\partial\Omega \cap
U$.
\end{itemize}
Then there do not exist any automorphism orbits in $\Omega$
accumulating at $p$.
\end{theorem}

Notice that this in particular implies
\begin{corollary} {\rm (Byun)} \sl The Kohn-Nirenberg Domain
$$
\Omega = \{ (z,w) \in \CC^2 : \Re w + |zw|^2 + |z|^8 +
\frac{15}7 |z|^2 \ \re z^6 < 0 \}
$$
has no automorphism orbit accumulating at the origin $(0,0)$.
\end{corollary}

\medskip
In the case that the automorphisms extend to the boundary, one can
say more about the Greene-Krantz conjecture.  Recent articles
\cite{LAN} by M. Landucci and \cite{BYG} by J. Byun and H.
Gaussier show the following:
\medskip

\begin{quote} \sl
If a domain satisfies Condition $R$ of Bell (that the Bergman
projection maps $C^\infty(\overline{\Omega})$ to itself---see
\cite{KRA3}), if a straight line segment of positive length lies
in its boundary, and if each point of the segment is convexifiable
and of minimum type (here infinite type is allowed), then none of
the point on the segment can be an orbit accumulation point. \rm
\end{quote}
\medskip

Despite such encouraging and supporting evidences, the most
general case of Greene-Krantz conjecture still awaits a solution.
At the same time we would like to pose the more restricted
problem:

\begin{Problem} \sl Let $\Omega$ be a bounded domain in $\CC^{n+1}$,
for some $n\ge1$, with $C^\infty$ smooth boundary. Then can one
show that an isolated infinite type boundary point cannot be an
orbit accumulation point?
\end{Problem}

\section{Asymptotic Behavior of Holomorphic Invariants}

\subsection{Boundary behavior of the Kobayashi and
Carath\'eodory metrics}

Methods of scaling have been used to study the boundary
asymptotics of invariant metrics on strongly pseudoconvex domains
in $\CC^{n+1}$. In what follows, if $\Omega_1$, $\Omega_2$ are
domains then we let $\Omega_2(\Omega_1)$ denote the collection of
holomorphic mappings from $\Omega_1$ to $\Omega_2$.  As usual, $D$
denotes the unit disc in $\CC$.  We begin by defining two
important invariant metrics in complex Finsler geometry (for
background and details, see \cite{KRA1}). These should be thought
of as generalizations of the Poincar\'{e} metric from the disc to
more general domains.
\begin{definition}  \rm
If $\O \ss \CC^{n+1}$ is open, then the {\it infinitesimal
Carath\'{e}odory} {\it metric} is given by $F_C: \O \times
\CC^{n+1} \ra \RR$ where
$$
 F_C(z,\xi) = \sup_{f \in B(\O)
                      \atop
                      f(z) = 0}
                                  |f_*(z) \xi|
              \equiv \sup_{f \in B(\O)
                      \atop
                      f(z) = 0}
      \left | \sum_{j=0}^{n} \frac{\partial f}{\partial z_j}(z)
              \cdot \xi_j \right | .
$$
\end{definition}
\begin{definition}  \rm
Let $\O \ss \CC^{n+1}$ be open.
Let $e_1 = (1,0,\dots,0) \in \CC^{n+1}.$
The infinitesimal form of the {\it Kobayashi/Royden metric}
is given by $F_K: \O \times \CC^{n+1} \ra \RR,$ where
$$
F_K(z,\xi)  \equiv  \inf\{|\a| ~:~ \exists f \in \O(B)\ \mbox{\rm
with} \ f(0) = z,
  \left ( f'(0)\right )(e_1) = \xi/\a \}.
$$
\end{definition}

Some alternative definitions are available.  For instance the
following invariant form is known:
\begin{eqnarray*}
F_K(z,\xi)
   & = & \inf \Big\{ \frac{|\xi|}{|(f'(0))(e_1)|} ~:~  \\
   & & \qquad f \in \O(B), (f'(0))(e_1)
   \mbox{ \rm is a constant multiple of} \ \xi \Big\}.
\end{eqnarray*}

Very little was known about these metrics---except for rather
abstract generalizations---until the seminal work of Ian Graham in
1975 (see \cite{GRA}).  His basic result is this:

\begin{theorem}[I. Graham] \sl
Let $\O \subset \! \subset \CC^{n+1}$ be a strongly pseudoconvex
domain with $C^2$ boundary. Fix $P \in \bO.$ Let $\xi \in
\CC^{n+1}$, and write $\xi = \xi_T + \xi_N,$ the decomposition of
$\xi$ into complex tangential and normal components relative to
the geometry at the point $P.$ Let $\rho$ be a defining function
for $\O$ normalized so that $|\nabla \rho(P)| = 1.$ Let
$\Gamma_\a(P)$ be a non-tangential approach region at $P$. If $F$
represents either the Carath\'{e}odory or Kobayashi/Royden metric
on $\O$ then
$$
 \lim_{\Gamma_\alpha (P) \ni z \ra P} d_\O(z) \cdot F(z,\xi)
 =  \frac{1}{2}|\xi_N| .
$$
\noindent Here $| \ \ |$ denotes Euclidean length and $d_\O(x)$ is
the distance of $x$ to the boundary of $\O.$
If $\xi = \xi_T$ is complex tangential, then we have
$$
\lim_{\Gamma_\a(P)\ni z \ra P} \sqrt{d_\O(z)} \cdot F(z,\xi) =
\frac{1}{2}{\mathcal L}(\xi,\xi) ,
$$
where
$$
{\mathcal L} (\xi,\xi) = \sum_{j, k=0}^n \frac{\partial^2
\rho}{\partial z_j \partial \bar z_k} \Big|_p \cdot \xi_j \xi_k
$$
for $\xi = (\xi_0, \ldots, \xi_n)$. Such $\mathcal L$ is called
the Levi form for the defining function $\rho$.
\end{theorem}

Graham's proof is based upon an intricate local analysis with
uniform estimates on the $\bar\partial$ operator.  Later, it has
turned out that the proof using the scaling methods are easier to
understand and gives finer analysis. Such subsequent analyses are
found in \cite{MA}, \cite{FU}, \cite{BOASY}, \cite{LEE},
\cite{ALA}, and others. Here we present S. Lee's refinement and
proof by the scaling method.
\medskip  \\

\begin{theorem} {\rm (S. Fu, D. Ma, S. Lee)} \sl
Let $\Omega$ be a bounded domain in $\CC^{n+1}$ with a $C^2$ smooth,
strongly pseudoconvex boundary. Let $F_\Omega (p,\xi)$ denote
either the Carath\'eodory or Kobayashi metric of $\Omega$ for $\xi
\in T_p \Omega = \CC^{n+1}$. For each $q \in \Omega$ sufficiently
close to the boundary of $\Omega$, choose a boundary point $p \in
\partial\Omega$ that is the nearest to $q$. Then it holds that
$$
\lim_{q \to \partial\Omega} %
\bigg( \Big( \frac{\|\xi_{N,p}\|}{2 d(q,\partial\Omega)}
 \Big)^2
+ %
\frac{L_{\partial\Omega, p} (\xi_{T,p},
\xi_{T,p}) }{d(q,\partial\Omega)} \bigg) %
\cdot F_\Omega (q,\xi)^{-2} %
= 1,
$$
where:
\begin{itemize}
\item[(1)] $d(q,\partial\Omega)$ is the distance between $q$ and
$\partial\Omega$,
\item[(2)] $\xi_{N, p}$ and $\xi_{T,p}$ denote the normal and the
tangential components to $\partial\Omega$ at $p$ of the vector
$\xi$, understood by a parallel translation as a vector at $p \in
\partial\Omega$, and %
\item[(3)] $L_{\partial\Omega, p}$ represents the normalized Levi
form of $\partial\Omega$ at $p$.
\end{itemize}
\end{theorem}

Notice that this result analyzes the asymptotic boundary behavior
of the Carath\'eodory and Kobayashi metric in all directions,
without restricting the trajectory of the point $q$ as it
approaches the boundary $\partial\Omega$.  Moreover, this shows
that the Carath\'eodory and Kobayashi metrics are asymptotically
Hermitian.
\medskip   \\

\noindent\\bf Sketch of the proof. \rm
Lee's proof proceeds following the scaling method.  Let $\Omega$
and $p \in \partial\Omega$ be as in the hypothesis of Theorem.
Following the original work of Graham \cite{GRA}, one first
observes that there exists an open neighborhood $U$ of $p$ for
which one has
$$
F_\Omega (q,\xi) \sim F_{\Omega \cap U} (q, \xi)
$$
as $q$ approaches $p$. Then, shrinking $U$ if necessary, apply the
scaling method to $\Omega \cap U$.  With the dilatation sequence
associated with $q$ above, which is the centering map followed by
a stretching map $\Lambda_j = L_j \circ A_j$, one sees that the
following hold:
\begin{itemize}
\item[(1)] With an appropriate fixed Cayley transform (a linear
fractional transformation) $\Phi$, the sequence of domains
$\Phi\circ \Lambda_j (\Omega \cap U)$ converges normally to the
open unit ball.
\item[(2)] $\Phi \circ \Lambda_j (q) = 0$ for every $j$.
\end{itemize}
Therefore, by the invariance and interior stability properties of
the Kobayashi and Carath\'eodory metrics, one sees that
\begin{eqnarray*}
F_{\Omega \cap U} (q, \xi)  =  F_{\Phi\circ\Lambda_j (\Omega
\cap U)} (0, d[\Phi\circ\Lambda_j]_q (\xi)) \\
 \sim  F_{\BB} (0, d[\Phi\circ\Lambda_j]_q (\xi)).
\end{eqnarray*}
Now a careful analysis of the last term yields the desired
conclusion. See \cite{LEE} for the detailed analysis. \hfill
$\Box$
\bigskip  \\

One easily sees in this method that strong pseudoconvexity is not
a restricting factor for this type of analysis.  Indeed, Lee in
the same paper showed how to analyze the boundary behavior of such
metrics in a domain with an exponentially flat, infinite type
boundary.

\subsection{Boundary Behavior of the Bergman Invariants}

The role of the dilatation sequence in the preceding section is
two-fold:
\begin{itemize}
\item[(1)] It turns the boundary limiting behavior problem into an
interior stability problem.
\item[(2)] The normal convergence limit of the sequence of sets
becomes generally simple.
\end{itemize}
Therefore study of boundary behavior of several holomorphic
invariants can be handled through scaling as long as they are
localizable and have interior stability.
\medskip     \\

Thus it is natural to reprove and refine the celebrated theorem of
Klembeck \cite{KLE} on the boundary behavior of Bergman curvature
in the strongly pseudoconvex domain by scaling.


The Bergman kernel, metric and curvatures have to be introduced.
For a bounded domain $\Omega$ in $\CC^{n+1}$, we consider the
square integrable holomorphic functions
$$
\mathcal{A}^2(\Omega) = \{ f :\Omega \to \CC : \hbox{
holomorphic}, \int_\Omega f d\mu < \infty\}
$$
where $d\mu$ denotes the standard Lebesgue measure for
$\CC^{n+1}$. It is obviously a linear subspace of the Lebesgue
space $L^2(\Omega)$. It follows by the Cauchy estimate that
$\mathcal{A}^2(\Omega)$ is a Hilbert space. Moreover, it is a
separable Hilbert space with respect to the standard $L^2$ inner
product.

Let $z$ be a point in $\Omega$. Consider the point evaluation map
$$
\Psi_z : \mathcal{A}^2(\Omega) \to \CC
$$
defined by $\Psi_z (f) = f(z)$.  The Cauchy estimates imply that
this is a bounded linear functional.  Consequently, the Riesz
representation theorem implies that there exists a holomorphic
function $k_z \in \mathcal{A}^2(\Omega)$ such that
$$
\Psi_z (f) = \int_\Omega f(\zeta) \overline{k_z(\zeta)} \
d\mu(\zeta)
$$
for every $f \in \mathcal{A}^2(\Omega)$. Let $K_\Omega (z,\zeta)
:= \overline{k_z (\zeta)}$ for $(z,\zeta) \in \Omega \times
\Omega$. It is known that the function satisfies the following
properties:
\begin{itemize}
\item[(\romannumeral1)] $K_\Omega(z,\zeta)$ is holomorphic in
$z=(z_0,\ldots,z_n)$ and conjugate holomorphic in $\zeta =
(\zeta_0,\ldots, \zeta_n)$.
\item[(\romannumeral2)] $K_\Omega (z,\zeta) = \overline{K_\Omega
(\zeta, z)}$.
\item[(\romannumeral3)] $\displaystyle{f(z) = \int_\Omega
K(z,\zeta)f(\zeta)\ d\mu(\zeta)}$ for every $f \in
\mathcal{A}^2(\Omega)$.
\end{itemize}

The function $K_\Omega$ is the celebrated {\it Bergman kernel
function} for $\Omega$. Bergman showed that the bilinear form
$$
\beta_p \equiv \sum_{j,k=0}^n \frac{\partial^2 \log K_\Omega
(z,z)}{\partial z_j \partial\bar z_k} \Big|_p \ dz_j \otimes d\bar
z_k
$$
defines a positive definite Hermitian form on the tangent space at
$p$ for $\Omega$. This is the Bergman metric. Following the
formalism in differential geometry this metric admits various
concepts of curvatures including the notion of holomorphic
sectional curvature.

Now we are ready to present the main theorem of this subsection.

\begin{theorem} \bf (Klembeck, Kim, Yu) \sl
Let $\Omega$ be a bounded domain in $\CC^{n+1}$ with a boundary
point $p \in \bO$.  Denote by $S(q,\xi)$ the holomorphic sectional
curvature of the Bergman metric of $\Omega$ at $q \in \Omega$ in
the direction $\xi \in T_q\Omega$. If there exists an open
neighborhood of $p$ in $\CC^{n+1}$ such that $\bO$ is
$\mathcal{C}^2$ strongly pseudoconvex at every point in $\bO \cap
U$, then
$$
\lim_{\Omega \ni q \to p} S(q,\xi) = -\frac4{n+2}.
$$
\end{theorem}


\noindent\bf Sketch of the proof. \rm It seems appropriate to
point out that the number $-4/(n+2)$ is actually the holomorphic
sectional curvature (from here on, we shall simply say {\it
holomorphic curvature}) of the Bergman metric at the origin of the
unit ball $B$ in $\CC^{n+1}$. Therefore Klembeck's theorem simply
says that the holomorphic curvature of a bounded domain is
asymptotically the holomorphic curvature of the unit ball at the
origin, as the reference point approaches a strongly pseudoconvex
boundary point. Originally, Klembeck proved the above stated
theorem with the stronger assumption that the domain has
$\mathcal{C}^\infty$ strongly pseudoconvex boundary. He needed
such a strong assumption because he used the celebrated asymptotic
expansion formula of C. Fefferman for the Bergman kernel function.
The improvement by Kim and Yu is in that they avoided using
Fefferman's formula by exploiting the convenience of the scaling
method for this type of problems, and that as a consequence they
could prove the theorem with $\mathcal{C}^2$ smoothness, the
optimal regularity assumption for strong pseudoconvexity.

The actual arguments by Kim and Yu proceed much along the line of
scaling methods demonstrated above. (In fact, this type of using
the scaling methods to study asymptotic behavior of holomorphic
invariants started with this problem in \cite{KI4}. See \cite{KIY}
for further developments.  Note also that Klembeck's result gives
another way to prove Wong-Rosay theorem.)

With the notation in the theorem, we let $q$ represent the general
element of a sequence of interior points of $\Omega$ approaching
$p$. We shall proceed much in the same way as in the scaling proof
of Graham's theorem in Section 6.2.

First one localizes the problem. Namely, for the sequence of
non-zero tangent vectors $\xi_q \in T_q \Omega$, one needs to
establish that for any open neighborhood $V$ of $q$, there exists
an open set $U$ with $p \in U \subset V$ such that the relation
$$
S_{\Omega} (q;\xi_q) \sim S_{\Omega \cap U} (q;\xi_q)
$$
holds as $q$ approaches $p$. This was done in \cite{KIY} in detail
(see also \cite{KIL}) using two classic results: the
representation of the holomorphic curvature of the Bergman metric
by the minimum integrals and the $L^2$ estimates of $\bar\partial$
operator by H\"ormander. It should be mentioned that the argument
by Kim-Yu does not use any regularity of the boundary but uses
only the existence of holomorphic peak functions at $p$ and the
pseudoconvexity of the domain $\Omega$ in consideration.

Then the next step is to change the holomorphic coordinates, by a
biholomorphism $\psi$ of $U$ onto an open ball $\widetilde U$ in
$\CC^n$ centered at the origin $0$ with an appropriate radius,
such that $\psi(p)=0$ and
\begin{eqnarray*}
\psi (\Omega \cap U) & = & \{z \in \widetilde U : \\
&& \quad \re z_0 > |z_0|^2 + \ldots + |z_n|^2  \\
&& \qquad + o(|z_0|^2 + \ldots + |z_n|^2) \} \, .
\end{eqnarray*}

Let $\widetilde q = \psi(q)$. Note that $\widetilde q$ approaches
the origin $0$ now. Then on $\psi(\Omega \cap U)$ with $\widetilde
q$, let us build the centering map $A_{\widetilde q}$, the
stretching map $L_{\widetilde q}$ and hence the scaling map
$\Lambda_{\widetilde q} \equiv L_{\widetilde q} \circ
A_{\widetilde q}$.

Now observe that the sequence of sets $\Lambda_{\widetilde q}
(\psi(\Omega \cap U))$ converges normally to the Siegel half space
$\mathcal{S}$ in $\CC^n$ defined by
$$
\re z_0 > |z_1|^2 + \ldots + |z_n|^2
$$
and $\Lambda_{\widetilde q} (\widetilde q) = (1,0,\ldots,0)$.

Let $\Phi$ be a standard Cayley linear fractional transformation
that maps $\mathcal{S}$ bihomorphically onto the unit ball $B$
such that $\Phi (1,0,\ldots, 0) = (0,\ldots, 0)$.

Since the holomorphic curvature of the Bergman metric is a
biholomorphism invariant, one immediately deduces that
\begin{eqnarray*}
S_\Omega (q, \xi_q) & \sim & S_{\Omega \cap U} (q, \xi_q) \\
& = & S_{\Phi \circ \Lambda_q (\Omega \cap U)} (\Phi \circ
\Lambda_q (q), d(\Phi \circ \Lambda_q)|_q (\xi_q)) \\
& \sim & S_{B} (0,d(\Phi \circ \Lambda_q)|_q (\xi_q)) \\
& = & - \frac{4}{n+2} \, ,
\end{eqnarray*}
which is the conclusion of the theorem. \hfill $\Box$

It should be obvious at this point that using the same method one
can study the boundary behavior of the Bergman metric as well as
the kernel itself in several cases including, but not limited to,
the strongly pseudoconvex domains. The interested reader may
consult the articles such as \cite{BOASY}, \cite{KIL} and the
references therein.

\subsection{Boundary Asymptotics of the Poisson Kernel}

In the recent paper \cite{KRA2}, Krantz uses a scaling method to
derive results on the boundary asymptotics of the Poisson kernel
on a bounded domain $\Omega \ss \RR^{n+1}$ with ${\mathcal C}^2$
boundary.  A typical result is

\begin{theorem}[Krantz] \sl
Let $\Omega \ss \RR^{n+1}$ be a bounded domain with $C^2$
boundary. Let $P: \Omega \times \partial \Omega \ra \RR^+$ be the
Poisson kernel for $\Omega$.  Let $\delta(x) = \delta_\Omega(x)$
denote the distance of $x$ to $\bO$.  Then there are constants
$c_1, c_2 \ge 0$ such that
$$
c_1 \cdot \frac{\delta(x)}{|x - y|^{n+1}} \leq P(x,y) \leq
      c_2 \cdot \frac{\delta(x)}{|x - y|^{n+1}} \, .    \eqno (\star)
$$
\end{theorem}

It is worth mentioning that the scaling sequence here is
isotropic, as the Laplace operator should be kept within the same
conformal class.

\section{Scaling in Infinite Dimensions}

In recent years the complex function theory of infinite
dimensions, particularly of Hilbert space, has received
considerable attention.  This new venue helps to put the classical
finite-dimensional situation into new perspective, and offers many
new challenges.  We offer here one example of a result that can be
proved by a scaling method, although it must be emphasized that
many aspects of the classical arguments of Bun Wong and Rosay as
well as a direct application of scaling methods fail in this
context, and many of them require new ideas:

\begin{theorem} \sl
Let $\Omega$ be a bounded, convex domain in a separable Hilbert
space ${\mathcal H}$.  Assume that $\Omega$ admits a boundary point
${\mathbf p} \in \bO$ at which
\begin{enumerate}
\item[{\bf (1)}] $\bO$ is ${\mathcal C}^2$ smooth and strongly
pseudoconvex in a neighborhood of ${\mathbf p}$, and \item[{\bf
(2)}] there exist ${\mathbf q} \in \O$ and $f_j \in \Aut(\O)$
($j=1,2,\ldots$) such that $f_j ({\mathbf q})$ converges to
${\mathbf p}$ in norm as $j \to \infty$.
\end{enumerate}
Then $\O$ is biholomorphic to the unit ball $\BB = \{z \in
{\mathcal H} : \|{\mathbf z}\| < 1 \}$.
\end{theorem}

We do not include the detailed construction of the scaling
sequence in infinite dimensions. The interested reader may consult
the article by Kim-Krantz \cite{KIK2} and also \cite{KMA}.

\section{Further Results}

It should be apparent by now that the scaling method is a powerful
tool for studying the asymptotic behavior of holomorphic
invariants as well as the non-compact orbits of a domain. At least
it should be obvious that the scaling method is related to many
other problems that have been open for some time and still need to
be studied.

\subsection{Linearization of Holomorphic Maps}

Linearization is one of the traditional problems in function
theory. The original problem is this: given a mapping from a
subset (containing the origin) of a Euclidean space into another
subset in the same space, can one find a new local coordinate
system at the origin for the Euclidean space under consideration
so that the original mapping becomes a linear map?  Let us
consider the case when the subset above is a domain, say $\Omega$,
the Euclidean space is the complex space $\CC^{n+1}$ and the map
$f:\Omega \to \CC^{n+1}$ is a holomorphic mapping preserving the
origin. The linearization problem now is asking whether there
exists a local biholomorphic mapping near the origin preserving
the origin, say $h$, such that $h \circ f \circ h^{-1}$ is the
restriction of a complex linear mapping from $\CC^{n+1}$ into
itself.

As studied earlier by Latt\`es, Poincar\'e, Dulac and others (see
\cite{STE} and the references therein), one may ask whether there
is a linear mapping $L:\CC^{n+1} \to \CC^{n+1}$ such that the
sequence $\psi_j \equiv L^{-j} \circ f^j$ converges uniformly on
compact subsets of $U$ to an injective holomorphic mapping, as $j$
tends to infinity. Here, $f^{j+1} = f^j \circ f$ for each $j=1,2,
\ldots$.  Let $f^0$ denote the identity map. Let $f^{-k} =
(f^{-1})^k$. \it If the answer to this last question is
affirmative\rm , let $\psi = \lim_{j\to\infty} L^{-j} \circ f^j$.
Then it follows that
\begin{eqnarray*}
\psi & = & \lim_{j\to\infty} L^{-j} \circ f^j \\
& = & L^{-1} \circ (\lim_{j\to\infty} L^{-j+1} \circ f^{j-1})
        \circ f \\
& = & L^{-1} \circ \psi \circ f.
\end{eqnarray*}
This shows that $\psi \circ f \circ \psi^{-1} = L$, and hence the
linearization problem is solved.

Therefore a key question is this:

\begin{quote} \sl
Does there exist a linear map $L$ such that the sequence $L^{-j}
\circ f^j$ forms a normal family?
\end{quote}

In lieu of a complete explication of this problem, we begin by
pointing out that one cannot help but notice the strong
resemblance of this problem to the scaling method. The conditions
for the convergence of the sequence $L^{-j} \circ f^j$ have been
studied extensively, from a time much earlier than the time of
Pinchuk's initial studies of complex scaling method. For instance,
in case the map $f$ is a contraction, in the sense that each
eigenvalue of the Jacobian matrix $df_0$ at the origin has modulus
less than 1, then the obstruction to the convergence is known---it
comes down to a resonance relation, that is in fact a collection
of finitely many algebraic equations between the eigenvalues of
$df_0$. A typical result is that if the eigenvalues of $df_0$ are
free of resonance relations, then the sequence $L^{-j} \circ f^j$
form a normal family in a neighborhood of $0$. Consequetly, the
map $f$ is linearizable.
\medskip

Focusing still on contractions, we consider the situation related
to domains with non-compact automorphism group; this of course is
one of the key subjects of the present article. Consider a real
hypersurface $M$ in $\CC^{n+1}$ passing through the origin $0$.
Let us assume that there exists a holomorphic mapping $f$ defined
in a neighborhood $U$ of $0$ mapping $U$ into $\CC^{n+1}$ such
that
\begin{itemize}
\item[(A)] $f(0)=0$,
\item[(B)] $f$ is a contraction,

and

\item[(C)] $f$ locally preserves $M$, i.e., $f(M \cap U) \subset
M$.
\end{itemize}
Now we ask whether $f$ can be linearized. As one readily observes,
we are asking whether the condition (C) can replace the
resonance-free condition.
\medskip

In some sense, in case $M$ is a real-analytic hypersurface that is
strongly pseudoconvex and not locally biholomorphic to part of the
sphere, then (C) does replace the resonance-free condition. (See
\cite{KRL}, \cite{KIS}, \cite{EZH}, \cite{STE}.) But then it is
not hard to deduce the following statement:

\begin{proposition} \sl
Let $M$ be a germ of real-analytic hypersurface in $\CC^{n+1}$
passing through the origin. Assume also that $M$ is strongly
pseudoconvex. If there exists a local biholomorphism $f$ of
$\CC^{n+1}$ defined in a neighbhorhood of $0$ preserving the
origin and mapping $M$ to $M$, then $M$ is biholomorphic to a germ
at $0$ of the hypersurface defined by
$$
\re z_0 = |z_1|^2 + \ldots + |z_n|^2.
$$
\end{proposition}

The proof follows by the linearization. We give a rough sketch
only. Expecting a contradiction, let us assume that $M$ is not
biholomorphic to a germ of the quadratic surface described above.
Then $f$ is linearizable, as mentioned above. Let $\re z_0 = \rho
(z_1, \ldots, z_n, \im z_1)$ denote the defining relation of $M$.
Let us assume, after a reduction, that $f$ itself is linear. Now
replace $f$ by its Jordan canonical form. Roughly speaking, the
components $(f_0, \ldots, f_n)$ of the map $f$ will satisfy
$$
\re f_0 = \rho (\im f_0, f_1, \ldots, f_n).
$$
Notice that $f$ is almost a diagonal linear map with the moduli of
the eigenvalues all less than one.  Give weights to the variables,
so that the weight for $z_0$ is 2, and the weight of $z_\ell$ is 1
for $\ell=1, \ldots, n$. Therefore an iteration of this process
will imply that all the monomial terms in the Tayler expansion of
$\rho$ with degree higher than 2 must vanish.  This is strongly
analogous to the scaling method described earlier in this paper.
But then the conclusion is that $M$ has to be defined by a
quadratic equation.  Since any strongly pseudoconvex hypersurface
defined by a quadratic equation with the prescribed weight has to
be linearly biholomorphic to part of sphere, we have arrived at a
contradiction.  \hfill $\Box$
\bigskip

This line of investigation gives rise to a number of interesting
questions. Notice in particular that the preceding proposition
gives a short and simple proof to part of following theorem (see
the discussion following for terminology):

\begin{theorem} {\rm (R. Schoen \cite{SCH})} \sl
Let $M$ be a $C^\infty$ strongly pseudoconvex CR manifold of
hypersurface type. If a point $p \in M$ admits a $CR$ automorphism
$f$ of $M$ such that $\lim_{j \to \infty} f^j (q) = p$ for every
point on $M$ except possibly one point, then $M$ is $CR$
equivalent to the sphere or the sphere minus one point.
\end{theorem}

The concept of abstract CR manifold, of hypersurface type or of
more general type, requires a precise introduction. We refer to
\cite{BOG} for details. Nevertheless, it is not so difficult to
picture what a CR manifold of hypersurface type should be.  First
consider a smooth real hypersurface $M$ (of real dimension $2n+1$)
in $\CC^{n+1}$. For each $p \in M$, the real (extrinsic) tangent
space $T_p M$ contains complex $n$ dimensional complex vector
subspace in it. This gives rise to a subbundle $\mathcal D$ of
$\mathcal{T}M$ with complex fibers, with real rank 1 transversal
distribution. Implementing this type of bundle with some
conditions called integrability in an abstract way formulates the
concept of CR manifold. See \cite{BOG}. Further, it is also
possible to define the Levi form from this abstract setting, and
hence the concept of strongly pseudoconvex CR manifold of
hypersurface type. Due to the famous embedding theorems of
strongly pseudoconvex CR manifolds of hypersurface type, except
for the case of dimension 3 or 5, the abstract strongly
pseudoconvex smooth CR manifolds of hypersurfaces type are locally
equivalent to smooth strongly pseudoconvex CR hypersurfaces in a
complex Euclidean space. (See \cite{KUR}, \cite{AKA}, \cite{WEB2}
and \cite{NIR}.)
\medskip

Continuing discussions from the above-stated theorem of Schoen, we
feel that it is natural at this juncture to pose the following:

\begin{Problem} \sl
Let $M$ be a $\calC^\infty$ smooth, strongly pseudoconvex $CR$
hypersurface of $\CC^{n+1}$ passing through the origin, not
locally biholomorphic to a sphere. If there is a local
biholomorphic mapping $f$ of $\CC^{n+1}$ preserving the origin and
the surface $M$, then show that $f$ is linearizable.
\end{Problem}

\subsection{$CR$ Hypersurfaces with Special $CR$ Automorphisms}

We shall continue with the discussion above. Focusing more on {\it
CR} hypersurfaces, it may be appropriate to point out that
non-degeneracy of the Levi form (which is often called Levi
non-degeneracy) is a more natural concept than strong
pseudoconvexity, at least for the {\it CR} hypersurfaces. Then the
representative model for such Levi non-degenerate hypersurfaces
should be hyperquadrics. On the other hand, it turns out that the
existence of a contracting holomorphic map that preserves the {\it
CR} hypersurface is a condition not in general restrictive enough
to conclude that the {\it CR} hypersurface with a contracting
holomorphic automorphism must be biholomorphic to a hyperquadric.
A correct condition has been found by Kim and Schmalz in
\cite{KIS}.

Call a local biholomorphic map $f$ at $0$ preserving a {\it CR}
hypersurface germ $M$ at $0$ in $\CC^{n+1}$ a {\it CR hyperbolic
automorphism} of $M$ if $df_0$ is expanding along the normal
direction to $M$ while $df_0$ is contracting along a complex
tangential direction. Then one has the following result.

\begin{theorem}[Kim/Schmalz] \sl
Let $(M,0)$ denote a germ of a real-analytic Levi non-degenerate
hypersurface at $0$ in $\CC^{n+1}$ ($n \ge 2$). If $M$ admits a
$CR$ hyperbolic automorphism, then $M$ is biholomorphic to a germ
of hyperquadric in $\CC^{n+1}$.
\end{theorem}

Again, it is expected that the $C^\infty$ version (or even $C^2$)
of this theorem should be true, but the more general result
remains open at this time.

\subsection{Existence of a One-parameter Family of Automorphisms}

Another interesting problem (communicated to us by Wu-Yi
Hsiang) that is closely related to the subject of this article
is as follows:

\begin{question} \sl
Let $\Omega$ be a bounded domain in $\CC^{n+1}$ with $\mathcal
C^\infty$ boundary. If there exists a sequence $\varphi_\nu$ of
automorphisms of $\Omega$ and a point $q \in \Omega$ such that the
orbit $\varphi_\nu (q)$ accumulates at a boundary point $p$ of
$\Omega$, then does the holomorphic automorphism group $\aut
\Omega$ contain a noncompact one-parameter subgroup?
\end{question}

Recall that the scaling method (when it applies) usually finds
that the given domain is biholomorphic to a domain represented by
an inequality of type
$$
\re z_0 > \psi (z_1, \ldots, z_n).
$$
Hence there is a 1-parameter translation automorphism subgroup
along the $\im z_0$ direction.  So the answer to the above
question is positive in several cases listed below in which the
scaling sequence always converges to a nice domain having such a
translation and/or dilation:
\begin{itemize}
\item $\Omega$ is convex and Kobayashi hyperbolic (but the
boundary need not be smooth (see \cite{FRA}, \cite{KI2}));
\item $\Omega$ is a subdomain of $\CC^2$ and bounded, or if
$\Omega \subset \CC^2$ is defined by a pluri-subharmonic defining
function, say $\rho$, with $\int_\Omega dd' \log \psi = +\infty$
(see \cite{BER}).
\end{itemize}

Since the convergence of scaling in the general case is quite
difficult, it seems reasonable to appeal to a different viewpoint.
In particular, we present a recent result that has exploited the
circle of ideas surrounding the linearization problem described
above.

\begin{theorem} [K.T. Kim/S.-y. Kim] \sl
Let $\Omega$ be a bounded pseudoconvex domain in $\CC^{n+1}$ with
a real-analytic boundary, and with no non-trivial complex analytic
variety. If there exists $f \in \aut (\overline{\Omega})$ such
that $f(p)=p$ for some $p \in \partial\Omega$ and such that $f$ is
a contraction at $p$, then $\dim_{\RR} \aut \Omega \ge 2$.
\end{theorem}

Before sketching the proof, we point out that the recent result
by K. Diederich and S. Pinchuk on the reflection principle
(\cite{DIP}) implies that every holomorphic automorphism of
such $\Omega$ extends holomorphically across its boundary (see
also \cite{HUAN}). Hence the assumption that $f$ should belong to
$\aut (\overline{\Omega})$ is not so restrictive in this case.

Now, we sketch the proof. The argument is computational, and is
rather intricate and involved. However, the nub of the proof is
that $f$ can be linearized to a diagonal matrix unless the the
boundary $\partial\Omega$ near $p$ is biholomorphic to a {\it CR}
hypersurface defined by a weighted homogeneous polynomial
function. Then it is shown in \cite{KI3} that the linearizability
again implies that the hypersurface has to be defined by a
weighted homogeneous polynomial.  Altogether, one can conclude
from this, using normal families argument, that the domain
$\Omega$ is biholomorphic to a domain in $\CC^{n+1}$ defined by a
weighted homogeneous polynomial.  But this latter domain admits
two non-compact 1-parameter families of automorphisms: dilation
and translation. Thus the conclusion of the theorem is obtained.
\hfill $\Box$
\bigskip

The study of automorphism groups is an instance of Felix Klein's
{\it Erlangen program}.  It provides an algebraic/geometric
invariant for distinguishing and comparing domains in complex
space or a complex manifold.  It is proving to be a powerful
tool in many aspects of geometric analysis and function theory.
Certainly the scaling method, which is an outgrowth of the theory
of automorphism groups, is finding use in subjects ranging
from partial differential equations to differential geometry
to complex variables.  We hope that this exposition will
spur further interest in this circle of ideas.

\vfill
\eject

\vspace{1in}

Kang-Tae Kim

Department of Mathematics

Pohang University of Science and Technology

Pohang 790-784 The Republic of Korea

{\tt (kimkt@postech.ac.kr)}

\vspace{.7in}

Steven G. Krantz

American Institute of Mathematics

360 Portage Avenue

Palo Alto, California 94306-2244 U. S. A.

{\tt (skrantz@aimath.org)}

\end{document}